\documentclass{article}

\usepackage{PRIMEarxiv}

\usepackage[utf8]{inputenc} % allow utf-8 input
\usepackage[T1]{fontenc}    % use 8-bit T1 fonts
\usepackage{hyperref}       % hyperlinks
\usepackage{url}            % simple URL typesetting
\usepackage{booktabs}       % professional-quality tables
\usepackage{amsfonts}       % blackboard math symbols
\usepackage{nicefrac}       % compact symbols for 1/2, etc.
\usepackage{microtype}      % microtypography
\usepackage{lipsum}
\usepackage{fancyhdr}       % header
\usepackage{graphicx}       % graphics

\usepackage{pifont}
\usepackage{amsfonts,amssymb,graphicx,listings,bm}
\usepackage{times}
\usepackage{subfig}
\usepackage{float}
\usepackage{verbatim}
\usepackage{amsfonts}
\usepackage{graphicx}
\usepackage{amscd}
\usepackage{multirow}
\usepackage{fancybox}
\usepackage{booktabs}
\usepackage{diagbox}
\usepackage{makecell}
   
\usepackage{amsmath}
\usepackage{amsthm}
\usepackage{amssymb}
\usepackage{amsfonts}
\usepackage{geometry}
\usepackage{pstricks-add}
\usepackage[all]{xy}
\usepackage{verbatim}
\usepackage{enumerate}
\usepackage{tikz}
\usepackage{algorithm}
\usepackage{algpseudocode}

\usetikzlibrary{shapes.geometric, arrows}
\tikzstyle{arrow} = [thick,->,>=stealth]

\newtheorem{theorem}{Theorem}

\graphicspath{{media/}}     % organize your images and other figures under media/ folder

\title{Multi-Scale Frequency-Enhanced Deep D-bar Method for Electrical Impedance Tomography }   

\author{
  Xiang Cao$^{a}$, Qiaoqiao Ding$^{a}$, Xiaoqun Zhang$^{a}$\\
  $^{a}$ School of Mathematical Sciences and Institute of Natural Sciences \\
  Shanghai Jiao Tong University \\
  Shanghai, 200240, CHINA\\
  \texttt{\url{{spawner,  dingqiaoqiao, xqzhang}@sjtu.edu.cn}} 
}

\begin{document}
\maketitle

\begin{abstract}
The regularized D-bar method is a popular method for solving Electrical Impedance Tomography (EIT) problems due to its efficiency and simplicity.  It utilizes the low-pass truncated scattering data in the non-linear Fourier domain to solve the associated D-bar integral equations, yielding a smooth conductivity approximation. However, the D-bar reconstruction often presents low contrast and resolution due to the absence of accurate high-frequency information and the ill-posedness of the problem. In this paper, we propose a deep learning-based supervised approach for real-time EIT reconstruction. Based on the D-bar method, we propose to utilize both  multi-scale frequency enhancement and spatial consistency for a high-image-quality reconstruction.  Additionally, we propose a fixed-point iteration for solving discrete D-bar systems on GPUs for fast computation. Numerical results are performed for both the continuum model and complete electrode model simulation on KIT4 and ACT4  datasets  to demonstrate notable improvements in absolute EIT imaging quality.  
\end{abstract}

% keywords can be removed

\section {\bf Introduction}
\subsection{Problem Setting}
\textit{Electrical Impedance Tomography} (EIT) stands as a remarkable imaging technique that has garnered significant attention in inverse problems with many applications. This non-invasive method harnesses the principles of electrical conductivity to explore the inner structure of objects and organisms, which has been used in fields like non-destructive medical monitoring \cite{adler2021electrical}, industrial detection \cite{york2001status} and geophysics \cite{abubakar2009inversion}. Specifically, EIT involves the use of distinct voltage patterns through electrodes placed on the object's surface and the subsequent measurement of resulting current distributions. The objective is to recover the underlying conductivity distribution of the subject. 

The mathematical model for the EIT problem is formulated as the well-known \textit{Calderón inverse problem} \cite{calderon2006inverse}. In our context, we suppose that $\sigma(z) \in L^{\infty}(\Omega)$ refers to the positive conductivity distribution defined on a bounded domain $\Omega$ in $\mathbb{R}^2$. For a Dirichlet boundary condition $f(z)$ applied on $\partial{\Omega}$, the resulting electrical potential field $u(z)$ satisfies the following second-order elliptical PDE:
\begin{equation}\hspace{2cm}
    \begin{array}{cc}
        \nabla \cdot[ \sigma(z) \nabla u(z)]=0, &\quad z \in \Omega, \\
         u(z)|_{\partial \Omega} = f(z), &\quad z \in \partial\Omega.
    \end{array}
\end{equation}
This EIT equation defines the so-called Dirichlet to Neumann (DtoN) map $\Lambda_{\sigma}$ on the boundary:
\begin{equation}\hspace{1cm}
    \Lambda_{\sigma}:\left.\left.u\right|_{\partial \Omega}\in H^{1/2}(\partial \Omega) \longrightarrow \sigma \frac{\partial u}{\partial \nu}\right|_{\partial \Omega}\in \widetilde{H}^{-1/2}(\partial \Omega), 
\end{equation}
where $\nu$ is the normal vector along  $\partial \Omega$, and $\widetilde{H}^{-1/2}(\partial \Omega)$ stands for the dual of the corresponding Sobolev space $H^{1/2}(\partial \Omega)$.
The \textit{Calderón inverse problem} aims to recover the conductivity distribution $\sigma(z)$ inside $\Omega$ from boundary measurements modeled by the DtoN map $\Lambda_\sigma$. In practice, only a finite-dimensional approximation of the DtoN map can be obtained for retrieving the conductivity distribution, which is highly sensitive to measurement noise and modeling errors \cite{mueller2012linear}. Given the inherent ill-posedness of the problem, we aim to leverage the data-driven method combined with the EIT mathematical model for robust and improved real-time EIT reconstruction.

\subsection{Related work}
The theoretical foundation of the Calderón inverse problem is centered on the existence and uniqueness of solutions. Sylvester and Uhlmann \cite{sylvester1987global} were the first to prove that for a smooth conductivity function $\sigma(z)$ satisfying $\sigma(z) > 0$ and $\sigma(z) \in C^2(\Omega)$, the conductivity distribution can be uniquely determined from the DtoN map based on the complex geometric optics (CGO) solutions. In $\mathbb{R}^2$, Nachman \cite{nachman1996global} extended this result and provided an explicit reconstruction algorithm for $\sigma(z)$ under the condition $\sigma(z) \in H^2(\Omega)$ established on the scattering theory. Further progress by Astala and Päivärinta \cite{astala2006calderon} relaxed the regularity requirement, showing uniqueness for $\sigma(z) \in L^\infty(\Omega)$ in $\mathbb{R}^2$ using the quasi-conformal mapping theory. 

Despite these theoretical advances, the Calderón inverse problem is severely ill-posed, as small perturbations in the boundary data can lead to large errors in the reconstruction of $\sigma(z)$. This ill-posedness stems from the high-frequency components of $\sigma(z)$ being poorly captured by the DtoN map, necessitating the use of regularization techniques \cite{zhou2015comparison,  jin2012reconstruction, chan2004level, chung2005electrical} to stabilize the inversion process. From an optimization standpoint, optimization-based methods necessitate high-precision numerical solvers \cite{soleimani2005improving, jin2017convergent} for the forward and adjoint models, as well as robust optimization strategies, such as the regularized Newton-type optimization methods \cite{islam2014electrical, tan2020electrical}, and subspace-based optimization methods \cite{chen2009subspace, wei2016two}. During the non-linear optimization process, the GMRES algorithm \cite{saad1986gmres, liu2018pyeit} is commonly used to efficiently solve the linearized sub-problem at each iteration. Note that one of the primary limitations of iterative methods lies in their computational demands, while non-iterative methods provide a reconstruction directly from measurement with more efficiency but lower resolution. The D-bar method \cite{mueller2012linear, isaacson2004reconstructions, liu2023deepeit} transforms the conductivity problem into a complex scattering problem via the D-bar integral equation, allowing for an explicit reconstruction of the conductivity $\sigma(z)$ from boundary data using a specially constructed scattering transform and its inverse. Other examples of the direct methods encompass the factorization method \cite{kirsch2007factorization, harris2022regularization}, the direct sampling method \cite{chow2014direct}, and the enclosure method \cite{ikehata2000reconstruction}, which also boast elegant mathematical foundations.

In addition to the classical algorithms in EIT, there has been a surge of interest in leveraging supervised learning techniques to solve the EIT problem. The CNN-based approach \cite{hu2019image} directly approximates the mapping from the noised DtoN matrix to the corresponding ground truth using a conventional CNN-based architecture for absolute imaging. The TSDL method \cite{ren2019two} features a two-stage architecture comprising pre-reconstruction and post-processing stages, where 
a CNN-based network with residual connections post-processes the one-step pre-reconstructed image from the measurements. Building upon the regularized D-bar reconstructions, the Deep D-bar method \cite{hamilton2018deep, hamilton2019beltrami} explores the U-Net architecture \cite{ronneberger2015u} as a post-processing network as well for the real-time reconstruction. Other supervised approaches, such as \cite{fan2020solving, guo2021construct, wei2019dominant, xu2023enhancing}, are also established on some classical methods with sophisticated network design for improving the reconstruction quality. Here, we should note that one of the learning-based approaches, the deep unrolling technique,
manages the high computational cost for solving the non-linear EIT problem, as it unrolls the governing forward and adjoint PDEs into the network multiple times. Besides, unsupervised learning methods based on physics-informed neural networks (PINN) \cite{bar2019unsupervised, raissi2019physics} and the deep image prior (DIP) \cite{liu2023deepeit} have been explored for solving electrical inverse problems. However, Bar \textit{et al.} \cite{bar2019unsupervised} only focused on conductivity inversion using the PINN, assuming full observation of the electrical potential $u(z)$, rather than utilizing boundary data. Liu \textit{et al.} \cite{liu2023deepeit} employed the DIP to parameterize the conductivity distribution $\sigma(z)$, which still requires evaluating the corresponding forward and adjoint PDEs for optimization. Given such limitations of unsupervised learning methods, we focus on exploring the use of supervised learning techniques for high-resolution and real-time EIT reconstruction. 

\subsection{Our Method}

\begin{figure}[!htp]
    \centering
    \includegraphics[width=15.5cm]{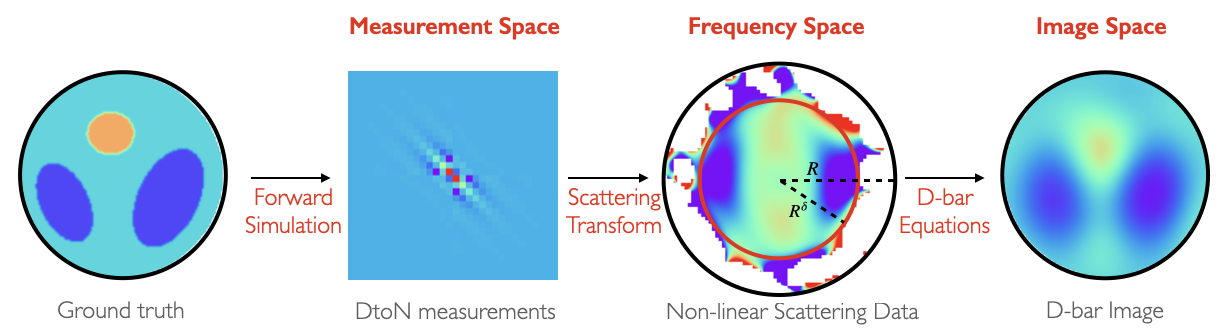}
    \caption[1]{From an insightful viewpoint, the D-bar regularization method belongs to direct reconstruction approaches for EIT.}
    \label{framework}
    \end{figure}

\textbf{Motivation.} Our motivation for investigating the learning-based D-bar method is twofold: first, its potential to achieve real-time computational efficiency, and second, its inherent regularization properties in the frequency domain. The D-bar method enables explicit reconstruction through the non-linear scattering transform and D-bar integral equations (as illustrated in Fig.~\ref{framework}). By eliminating the need to solve both forward and adjoint PDEs during optimization, it offers the potential for real-time reconstruction. Additionally, owing to its capability to partition the non-linear frequency domain into stable and unstable parts \cite{mueller2012linear}, the D-bar method employs low-pass filtering as a regularization technique in the frequency domain. However, the regularized D-bar reconstruction tends to poorly capture the high-frequency components of $\sigma(z)$, motivating the use of data-driven approaches to recover the high-frequency information and enhance the reconstruction process.

\textbf{Main Method.} Based on the considerations outlined above, we propose a multi-scale frequency-enhanced deep D-bar method for real-time EIT reconstruction. First, a GPU-based fixed-point iteration algorithm is implemented to solve the D-bar integral equations, enabling the real-time production of low-pass truncated D-bar images and frequency-enhanced D-bar images. Second, a cascaded multi-scale U-Net architecture is utilized for supervised post-processing. Specifically, the first U-Net predicts multi-scale frequency-enhanced D-bar images from low-pass truncated D-bar inputs, while the second U-Net calibrates these frequency-enhanced D-bar predictions for high-resolution EIT reconstruction. Additionally,  Numerical experiments are conducted using simulated KIT4 and ACT4 datasets based on the continuum model. For simulating real-world applications, the complete electrode model is considered for testing the model performance on KIT4 and ACT4 phantoms.

\textbf{Organization.} The remainder of the paper is structured as follows. Section~\ref{2} provides a concise overview of the regularized D-bar method. Building on the principles of the D-bar method, Section~\ref{3} introduces our proposed multi-scale frequency-enhanced deep D-bar method for real-time EIT reconstruction. Sections~\ref{4} through~\ref{6} present detailed experiments, numerical results, and ablation studies. Finally, Section~\ref{7} summarizes the key findings and conclusions of the study.

\section{\bf Regularized D-bar Method}\label{2} 

The traditional regularized D-bar method \cite{mueller2012linear} is summarized in the following three main steps: 
\begin{itemize}
    \item For each complex scattering variable $k\in\mathbb{C}$, we evaluate the so-called \textit{Complex Geometrical Optics} (CGO) solutions $\psi(z, k)$ on the boundary $\partial \Omega$ from the DtoN map $\Lambda_\sigma$:
    \begin{equation}\hspace{1cm}\label{CGO}
        \left.\psi(\cdot,k)\right|_{\partial \Omega}=\left.e^{i k \cdot}\right|_{\partial \Omega}-\delta_k * \left(\Lambda_\sigma-\Lambda_1\right) \psi(\cdot, k),
    \end{equation}
    where $\Lambda_1$ is the DtoN map for the background conductivity $\sigma(z) = 1$, and $\delta_k$ is defined as the Faddeev Green's function for the Laplacian $\Delta$. 
    
    \item Integrate the CGO solutions $\psi(z,k)$ along the boundary $\partial \Omega$ to obtain the associated scattering data $\mathbf{t}(k)$ from 
    \begin{equation}\hspace{1.6cm}\label{scattering}
        \mathbf{t}(k) = \int_{\partial \Omega} e^{i \bar{k} \bar{z}}\left(\Lambda_\sigma-\Lambda_1\right) \psi(z, k) d s.
    \end{equation}
    Due to the measurement errors in the DtoN map $\Lambda_{\sigma}$ and the exponential factor $ e^{i \bar{k} \bar{z}}$ in Eq.~(\ref{scattering}), $\mathbf{t}(k)$ becomes numerically unstable in the high-frequency region, necessitating the low-pass filtering as regularization for stability.
    
    \item For each $z \in \Omega$, we solve the regularized D-bar solutions $m_{r}(z, k)$ from the following integral equations:
    \begin{equation}\hspace{-0.6cm}\label{D-bar-equation}
        m_{r}(z, k)=1+\frac{1}{(2 \pi)^{2}} \int_{|k^{\prime}|<r} \frac{\mathbf{t}(k^{\prime})}{\left(k-k^{\prime}\right) \bar{k}^{\prime}}  e^{-i(k^{\prime}z + \bar{k}^{\prime}\bar{z})} \overline{m_{r}\left(z, k^{\prime}\right)} d k^{\prime},
    \end{equation}
    where $r$ is selected as the truncation radius such that the unstable part of $\mathbf{t}(k^{\prime})$ for $|k^{\prime}| > r$ is filtered. Finally, the regularized D-bar reconstruction $\sigma_{r}(z)$ comes from 
    \begin{equation}\hspace{3.2cm}\label{recon}
        \sigma_{r}(z)=m_{r}(z,0)^2.
    \end{equation}
\end{itemize}
As shown in Fig.~\ref{framework}, the reconstruction process $\Lambda_\sigma \to \mathbf{t}(k) \to \sigma_{r}(z)$ demonstrates that the scattering data $\mathbf{t}(k)$ is essential for linking the measurements to the conductivity distribution. Using the Alessandrini's identity \cite{nachman1996global}, the scattering data $\mathbf{t}(k)$ can also be directly obtained from the conductivity $\sigma(z)$ and the CGO solutions $\psi(z,k)$ through the so-called \textit{scattering transform}: 
\begin{equation}\hspace{2.5cm}\label{scattering_sigma}
	\mathbf{t}(k) = \int_{\Omega} q(z) e^{i \bar{k} \bar{z}} \psi(z,k) d z,
\end{equation}
where $q(z) := \Delta \sqrt{\sigma(z)}/\sigma(z)$. The corresponding first-order approximation $\tilde{\mathbf{t}}(k)$ is given by 
\begin{equation}\label{scattering_approx}\hspace{2.5cm}
     \tilde{\mathbf{t}}(k) := \int_{\Omega} q(z) e^{i (\bar{k} \bar{z} + k z)} dz,
\end{equation}
and is commonly used to approximate $\mathbf{t}(k)$ in the high-frequency region due to the asymptotic behavior of the CGO solutions $\psi(z,k) \to e^{i k z}$ as $|k| \to \infty$.

\textbf{Deep D-bar Method.} To improve the reconstruction quality, the learning-based spatial domain post-processing method was introduced in \cite{hamilton2018deep, hamilton2019beltrami}. In this approach, the U-Net \cite{ronneberger2015u} is directly applied to the regularized D-bar reconstruction $\sigma_{r}(z)$. The well-established U-Net $f_{\Theta}(\cdot)$ is employed to model the relationship between the regularized D-bar reconstruction $\sigma_{r}(z)$ and the ground truth $\sigma(z)$ by minimizing the $L_2$ loss function 
\begin{equation}\hspace{2.3cm}\label{U-Net-loss}
      \mathcal{L}(\Theta) := \sum\limits_{n = 1}^{N} \|f_{ \Theta}(\sigma^{(n)}_{r}) - \sigma^{(n)}\|^2_{2},
\end{equation}
where $\Theta$ is the learnable parameter of the neural network and $n$ is index of $N$ training samples. Note that selecting the truncation radius $r$ for the network input $\sigma_{r}(z)$ holds immense significance in determining the generalization ability of the post-processing model. An ill-suitable choice of $r$ can lead to low-contrast D-bar reconstructions or deceptive artifacts. This is primarily because it predominantly preserves low-frequency information, while high-frequency data is more susceptible to amplified noise. Building upon the deep D-bar framework, we further investigate a more sophisticated learning framework design to improve the model performance.

\section{\textbf{Proposed Method}} \label{3}
%A uniform truncation radius $r$ in  Eq.~(\ref{D-bar-equation}) ensures contrast-consistent D-bar images $\sigma_r(z)$ across different conductivity samples $\sigma(z)$. 

\begin{figure}[!htp]
   \centering
   \includegraphics[width=15cm]{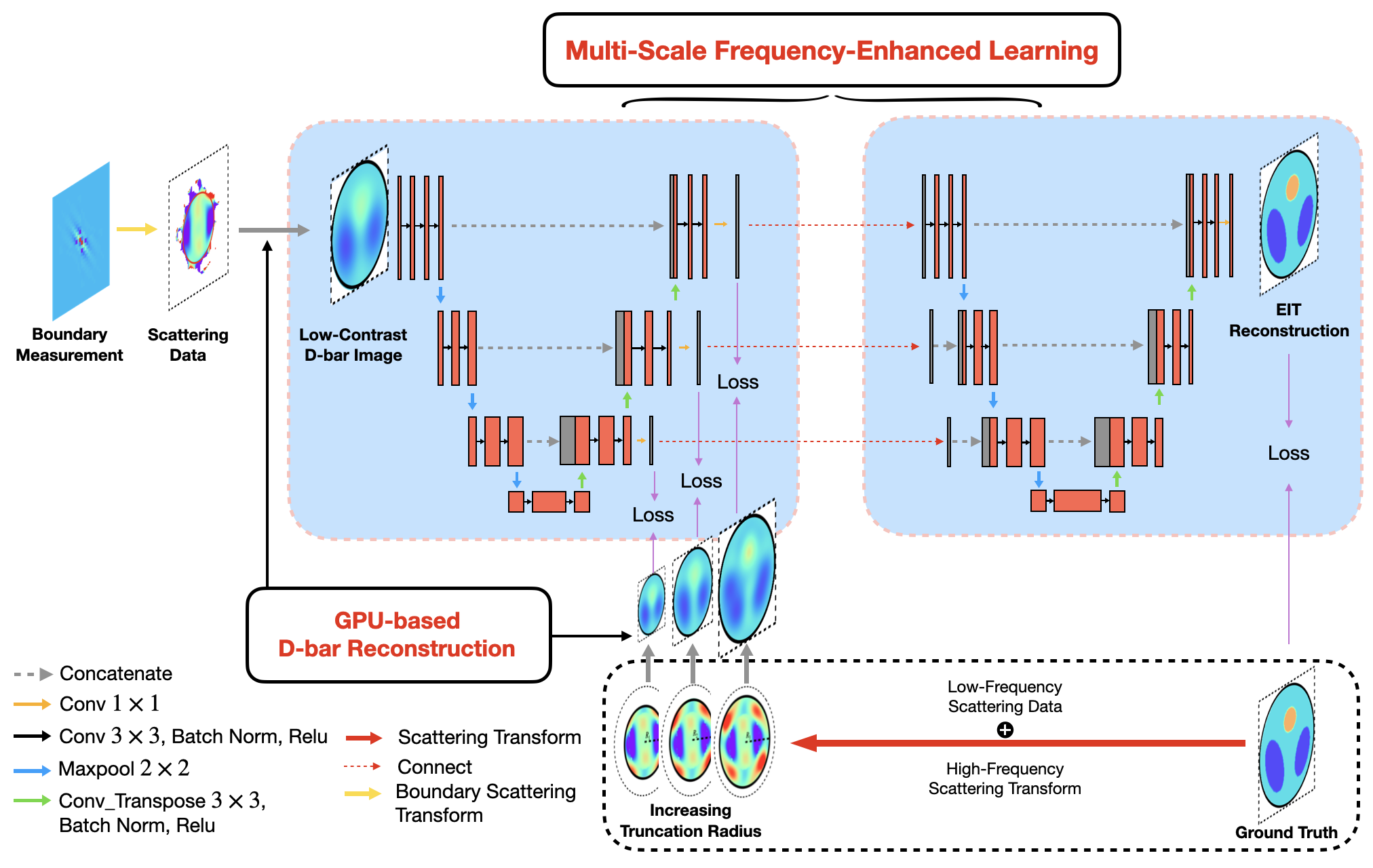} 
   \caption[]{The multi-scale frequency-enhanced learning framework consists of two cascaded blocks for EIT reconstruction. Additionally, the GPU-based D-bar reconstruction is introduced for real-time computation.}
  \label{DDDD}
\end{figure}

To improve the performance of the learning-based post-processing model, a larger radius $r$ is preferred to enhance the contrast of D-bar images, as it incorporates more high-frequency scattering information. However, due to the measurement noise and the inherent numerical instability, the scattering data $\mathbf{t}(k)$ in Eq.~(\ref{scattering}) is unavailable for $|k| > r$. This limitation necessitates a data-driven approach to recover the missing high-frequency scattering data to produce high-contrast D-bar images $\{\sigma_{r_i}(z)\}_{i=1}^{I}$ for $r_i > r$. In this work, we choose to directly learn to predict these images in the image domain, rather than learn to complete the missing high-frequency scattering data in the frequency domain. In this section, we will introduce the proposed multi-scale frequency-enhanced learning framework, as illustrated in Fig.~\ref{DDDD}, characterized by both multi-scale frequency and image enhancement with two cascaded neural network blocks. 

\textbf{Multi-Scale Frequency-Enhanced Learning.} In our cascaded learning framework, we consider the frequency-enhancement part for the high-contrast D-bar image, followed by the image-calibration part for the ultimate EIT reconstruction. Here, we employ the commonly used U-Net architecture to construct the frequency-enhancement network $f_{\Theta_1}(\cdot)$ and the image-calibration network $f_{\Theta_2}(\cdot)$. For clarity, we first define the following components that represent the U-Net architecture:
\begin{itemize}
    \item \textit{Encoder Blocks}: Each encoder block $\mathbf{E}_i$ maps the higher-level feature $X_{i-1}$ to the lower-level feature $X_{i}$, 
    $$ X_{i} = \mathbf{E}_i(\downarrow X_{i-1}),$$
    where $\downarrow$ is the downsampling operation. 
    \item \textit{Bottleneck}: The bottleneck $\mathbf{B}$ is the deepest layer, 
    $$ Y_{\tilde{I}} = \mathbf{B}(X_{\tilde{I}}), $$
    where $\tilde{I}$ is the number of encoder blocks.
    \item \textit{Decoder Blocks}: Each decoder block $\mathbf{D}_i$ maps the lower-level feature $Y_{i+1}$ to the higher-level feature $Y_{i}$, 
    $$ Y_{i} = \mathbf{D}_i(\uparrow Y_{i+1} \oplus X_{i}), $$
    where $\uparrow$ is the upsampling operation, and $\oplus$ represents the concatenation. 
    \end{itemize}
Note that the target frequency-enhanced D-bar images $\{\sigma_{r_i}(z)\}_{i=1}^{I}$ provide multi-scale scattering information from multiple truncation radii. To integrate such information into the cascaded learning framework, we design a multi-scale skip connection between two cascaded learning modules. Specifically, we customize the architecture of the frequency-enhancement network $f_{\Theta_1}(\cdot)$ by introducing an additional output layer at each decoder level. Specifically, we introduce the output layer $\mathbf{O}_i(\cdot)$ to the decoder output feature $Y_i$ and define
\begin{equation}\hspace{3.5cm}
    f^{(i)}_{\Theta_1}(\sigma_{r}) = \mathbf{O}_i(Y_i)
\end{equation}
as the $i$-level output of the first module. Meanwhile, the image-calibration network $f_{\Theta_2}(\cdot)$ is designed to encode the predicted multi-scale D-bar images at each encoder level. To be specific, we introduce the input layer $\mathbf{I}_i(\cdot)$ and define
\begin{equation}\hspace{2.5cm}
    X_{i} = \mathbf{E}_i\Big(\downarrow X_{i-1} \oplus \mathbf{I}_i\big(f^{(i)}_{\Theta_1}(\sigma_{r})\big)\Big) 
\end{equation}
as the $i$-level input of the second module. Finally, the combined supervised loss function is defined as 
\begin{equation}\hspace{-1.5cm}
 \mathcal{L}(\Theta_1, \Theta_2)  =  \sum\limits_{n = 1}^{N} \left(\Big\| f_{\Theta_2}\big[ f_{\Theta_1}(\sigma^{(n)}_{r})\big] - \sigma^{(n)} \Big\|^2_{2} + \alpha \sum\limits_{i = 1}^{I}\left\|f^{(i)}_{\Theta_1}(\sigma^{(n)}_{r}) - \sigma^{(n)}_{r_i} \right\|^2_{2} \right). 
\end{equation}
where $\alpha$ is a constant weight, and $N$ refers to the total number of training samples. % We further demonstrate the effectiveness of the multi-scale architecture compared to the single-scale architecture, as discussed in Sections~\ref{5} and~\ref{6}. 

\textbf{GPU-Based D-bar Reconstruction.} To generate the target high-contrast D-bar image $\sigma_{r_i}(z)$ from the ground truth $\sigma(z)$ for supervised learning, we evaluate the scattering transform $\mathbf{t}(k)$ in Eq.~(\ref{scattering_sigma}) and the D-bar equations in Eq.~(\ref{D-bar-equation}) to obtain the high-contrast D-bar image $\sigma_{r_i}(z)$ for a chosen truncation radius $r_i > r$. However, this  process is computationally expensive due to the need to solve the CGO solutions $\psi(z,k)$ for each $k$, and the regularized D-bar solutions $m_{r_i} (z,k)$ for each $z$. To achieve the real-time D-bar reconstruction, we propose the following two remedies for acceleration:
\begin{itemize}
    \item In the high-frequency part, we approximate $\mathbf{t}(k)$ for $r < |k| < r_i$ using the first-order approximation $\tilde{\mathbf{t}}(k)$ from Eq.~(\ref{scattering_approx}), obviating the need to solve numerical CGO solutions $\psi(z, k)$. In the low-frequency part, we directly use the pre-computed $\mathbf{t}(k)$ for $|k| < r$ from Eq.~(\ref{scattering}) during the low-pass D-bar reconstruction. 
    
    \item The fixed-point iteration is applied to solve the D-bar integral equations in parallel at different $z$-points on the GPU. See more details in the numerical implementations and the convergence proof of the fixed-point iteration in \ref{proof}.
\end{itemize}
Given the target high-contrast D-bar images $\{\sigma_{r_i}(z)\}_{i=1}^{I}$ generated from multiple truncation radii $\{r_i\}_{i=1}^{I}$ (selected as a monotonic decreasing sequence), we leverage such frequency-enhanced D-bar images as guidance in the aforementioned learning framework.

\section{\bf{Numerical Implementations}}\label{4}
In this section, we present our numerical implementations for the proposed learning-based method. The regularized D-bar method has been implemented in \cite{mueller2012linear} (see the provided \href{https://wiki.helsinki.fi/display/mathstatHenkilokunta/Inverse+Problems+Book+Page}{Matlab code} for computational resources). In this work, we utilize parts of the implementation, including the generation of synthetic conductivity phantoms, the simulation of the Neumann-to-Dirichlet matrix $\mathbf{L}_{\sigma}$, and the computation of low-pass truncated scattering data $\mathbf{t}(k)$ for $|k| < r$.

\subsection{\textbf{Boundary Measurement Models for Dirichlet-to-Neumann Map Simulation}}\label{4.1}
Two types of boundary measurement models are considered in the numerical experiments, \textit{the continuum model} and \textit{the complete electrode model}. The continuum model is applied for training and evaluating the framework using synthetic data, whereas the complete electrode model is used to test the learning-based model's performance for real-world applications. 

Specifically, we first compute the Neumann-to-Dirichlet (NtoD) map $\mathcal{R}_{\sigma}$ using $2P$ basis functions $\varphi_p(\theta)$ with $-P \le p \le P, p\ne 0$. For the continuum model, $\varphi_p(\theta)$ is chosen as the sinusoidal functions:
\begin{equation}\hspace{2cm}
     	\varphi_p(\theta)=\left \{
    \begin{array}{ll}
    	\frac{1}{\sqrt{2\pi}} \sin (p \theta) &\textrm{if  } p<0,\\
    	\frac{1}{\sqrt{2\pi}} \cos (p \theta)  &\textrm{if  } p>0.\\
\end{array}
\right .\nonumber
\end{equation}
For the complete electrode model, there are $2P$ equidistributed electrodes attached on the boundary, and $I^{l}_{p}$ is the current applied on the $l$-th electrode $e^l$ as:
\begin{equation}\hspace{2.2cm}
     I^{l}_p=\left \{
    \begin{array}{ll}
        \frac{1}{\sqrt{2\pi}} \sin(p\theta_l), &\textrm{if  } p<0,\\
        \frac{1}{\sqrt{2\pi}} \cos(p\theta_l), &\textrm{if  } p>0,\\
    \end{array}
    \right .\nonumber
\end{equation}
where $\theta_l$ is the angle of the midpoint of the corresponding electrode $e^l$. Then, $\varphi_p(\theta)$ is chosen as the piecewise constant functions: 
\begin{equation}\hspace{2.4cm}
     	\varphi_p(\theta)=\left \{
    \begin{array}{ll}
    	I^{l}_{p} &\textrm{if  } \theta \in e^{l},\\
    	0  &\textrm{if  }  \theta \notin \bigcup_{l=1}^{2P} e^{l}.\\
\end{array}
\right .\nonumber
\end{equation}

Given the boundary basis functions $\varphi_p(\theta)$, the matrix representation $\mathbf{R}_{\sigma} = [\mathbf{R}_{\sigma}]_{2P\times 2P}
 $ of the NtoD map $\mathcal{R}_{\sigma}$ can be approximated by 
\begin{equation}\hspace{1cm}
     [\mathbf{R}_{\sigma}]_{p,q} = \langle \mathcal{R}_{\sigma} \varphi_q, \varphi_p \rangle:= \int_{0}^{2\pi} (\mathcal{R}_{\sigma} \varphi_q)(\theta)\varphi_p(\theta) d\theta 
 \end{equation}
where  $\mathcal{R}_{\sigma} \varphi_p = \left.u_p\right|_{\partial \Omega}$ \text{with} $\nabla \cdot \sigma \nabla u_p=0 $ in $\Omega$, and the Neumann boundary condition $\left.\sigma\left(\frac{\partial u_p} {\partial \nu}\right)\right|_{\partial \Omega}= \varphi_p$ on $\partial \Omega$. The finite element method (FEM) is used to solve the EIT equations for each boundary condition. Furthermore, we consider the noised boundary voltage data $\mathcal{R}^{\delta}_{\sigma} \varphi_p$ defined as 
\begin{equation}\label{noise}\hspace{2cm}
     \mathcal{R}_{\sigma}^{\delta} \varphi_p := \mathcal{R}_{\sigma} \varphi_p +\delta \mathcal{N}_p\left\|\mathcal{R}_{\sigma} \varphi_p\right\|_{L^{\infty}}, 
\end{equation}
where $\delta$ denotes the simulated noise level, and $\mathcal{N}_{-P}, ...,\mathcal{N}_{P}$ are $2P$ independent noise samples from the standard Gaussian distribution. The matrix representation $\mathbf{L}_{\sigma}^{\delta}$ of the DtoN map $\Lambda^{\delta}_{\sigma}$
is generated by inverting $\mathbf{R}^{\delta}_{\sigma}$ and adding a zero column and row to the middle of it as described in \cite{mueller2012linear}. 
Subsequently, the numerical CGO solutions in Eq.(\ref{CGO}) are further expanded based on the basis functions $\varphi_p$, enabling the formulation of the corresponding linear systems. 
 
In this paper, $32$ different current patterns ($P = 16$) are applied to the boundary $\partial{\Omega}$ and the simulated noise levels $\delta = 0\%, 0.1\%, 0.75\%$ are considered for the noised measurement matrices $[\mathbf{L}_{\sigma}^{\delta}]_{33 \times 33}$. The corresponding truncated radius for the scattering data $\mathbf{t}(k)$ is chosen as $r = 6, 5, 4$ to guarantee the numerical stability in the $k$-space. 

\subsection{\textbf{Evaluation of the Scattering Transform and the D-bar Equations}}\label{Num_imp_D_bar_1}

The numerical implementation of the regularized D-bar method is referred to \cite{mueller2012linear}. Here, we summarize the numerical settings of evaluating the scattering transform and solving the D-bar equations. 

\textbf{Computational Grids.} The discrete points for the scattering transform $\mathbf{t}(k)$ are generated using $h = 0.2$ equidistant grids in the $k$-space. Given the numerical CGO solutions to $\psi(z, k)$ from Eq.~(\ref{CGO}), the scattering data $\mathbf{t}(k)$ for $ |k| < r $ is obtained from the boundary integral in Eq.~(\ref{scattering}). To evaluate the first-order approximation $\tilde{\mathbf{t}}(k)$ of the scattering transform for $r < |k| < r_1$, $269 \times 269$ equidistant $z$-grids in the square $[-2.1, 2.1)^2$ is generated to discretize $q(z)$ in Eq.~(\ref{scattering_sigma}). To avoid singularities, $q(z)$ is evaluated from the discontinuous conductivities $\sigma(z)$ with a slight Gaussian smoothing. 

Meanwhile, the D-bar reconstruction requires the computational grids generated on $512 \times 512$ equidistant $k$-grids in the square $[-2.1r_1, 2.1r_1)^2$, where $r_1$ is the chosen increased truncation radius. Then, the scattering data $\mathbf{t}(k)$ on the $k$-grids is obtained from the bilinear interpolation. Here, we choose three increased truncation radii $r_1, r_2, r_3 = 8, 7, 6$ to reconstruct the multi-scale frequency-enhanced D-bar images on $128 \times 128$, $64 \times 64$ and $32 \times 32$ equidistant $z$-grids in $[-1,1)^2$,  respectively. 

\textbf{Parallel Computation for D-bar Integral Equations.} To solve the numerical D-bar equations in real-time, we implement the fixed-point iteration to solve Eq.~(\ref{D-bar-equation}) leveraging parallel computation on the GPU. Specifically, 
starting with the initial condition $m^{(0)}_r(z,k) = 1$, we iterate the following linear operation:
\begin{equation}\label{D-bar-equation-iter}\hspace{1.5cm}
    m^{(s+1)}_r\left(z, k\right) = 1 + \left[\frac{1}{k}\right] \mathbin{\ast} \mathcal{T}_r\left[m^{(s)}_r\left(z, k\right) \right],
\end{equation}
where $\mathcal{T}_r[\cdot](z,k) = \frac{\mathbf{t}\left(k\right)\mathbb{1}_{|k|<r}(k)}{(2 \pi)^{2}\bar{k}}e^{-i(kz + \bar{k}\bar{z})} \overline{[\cdot]}$, and the convolution operation $\mathbin{\ast}$ can be efficiently evaluated via the Fast Fourier Transform. Compared to the GMRES algorithm used in \cite{mueller2012linear}, the fixed-point iteration method is more suitable for parallelizing the computation of D-bar solutions, $m_r(z,k)$, at different $z$-points on the GPU through concatenation. The numerical experiments in Section~\ref{6} demonstrate that, within just 5 iterations, the fixed-point iteration method can provide a reasonably accurate estimation of frequency-enhanced D-bar images. This approach offers substantial speedup compared to the CPU-based implementation using the GMRES algorithm in \cite{mueller2012linear}. 

\subsection{Network Implementation}

\textbf{Network Architecture.} The backbone U-Net architecture consists of multiple components arranged in a hierarchical structure to effectively handle multi-scale features. The input image, with one channel, is first projected into a feature space using a convolutional layer with a kernel size of $3$ and a hidden channel size of $32$. The encoder (downsampling path) comprises five resolution levels, where the number of channels increases according to the channel multipliers $[1,2,2,2,4]$. Each level includes three DownBlocks, each containing a ResidualBlock with two convolutional layers (kernel size $1,3$), shortcut connections ($1 \times 1$ convolution), GroupNorm, and GELU activation. Between levels, Downsample layers with $3 \times 3$ convolutions and stride $2$ reduce the spatial resolution by half. At the bottleneck, a MiddleBlock is used, composed of a ResidualBlock, and another ResidualBlock. The decoder (upsampling path) mirrors the encoder, starting with the highest resolution features and gradually increasing spatial resolution using Upsample layers, which apply transposed convolutions (kernel size $4$, stride $2$). At each resolution, UpBlocks are applied, each consisting of a concatenation of skip connections from the encoder and an UpBlock with ResidualBlock. The output of the final resolution is passed through a normalization layer and a final convolutional layer (kernel size $1$) to produce the output image with one channel.

\textbf{Training Configuration.} For the training configuration, an Adam optimizer is employed with a learning rate of $0.001$ and a weight decay of $10^{-6}$, while a step-wise scheduler adjusts the learning rate per epoch with a decay factor of $0.99$. The model is trained using the MSE loss function over $3000$ epochs, ensuring sufficient iterations for convergence. Multi-GPU data-parallel training accelerates the computation and accommodates a batch size of $128$ for training and $32$ for validation, enabling efficient utilization of available GPU memory while maintaining model performance.

\section{\textbf{Numerical Results}}\label{5}

\subsection{\textbf{Phantom Simulation for KIT4 and ACT4 Datasets}}
The KIT4 dataset is a collection of simulated phantoms consisting of several discontinuous inclusions with some regular shapes. We first create randomly selected anomaly shapes and locations in the interior. The inclusions are not allowed to superimpose, and each one is assigned an equal probability to be 'conductive' or 'resistive' for the background. Specifically, the values of the 'conductive' parts and the 'resistive' parts are sampled from the uniform distributions $[1.5, 2.5]$ and $[0.3, 0.7]$, respectively. The background's value is always set to be $1$, which could be naturally extended to the outside of the circle. A total of \textbf{3280 training samples and 820 testing samples} are simulated for each noise level $\delta$, and pairs of the scattering data $\mathbf{t}(k)$ and the ground truth $\sigma(z)$ are collected for training the neural network. 

\begin{figure}[!htp]
  \centering
  \includegraphics[width=15.5cm]{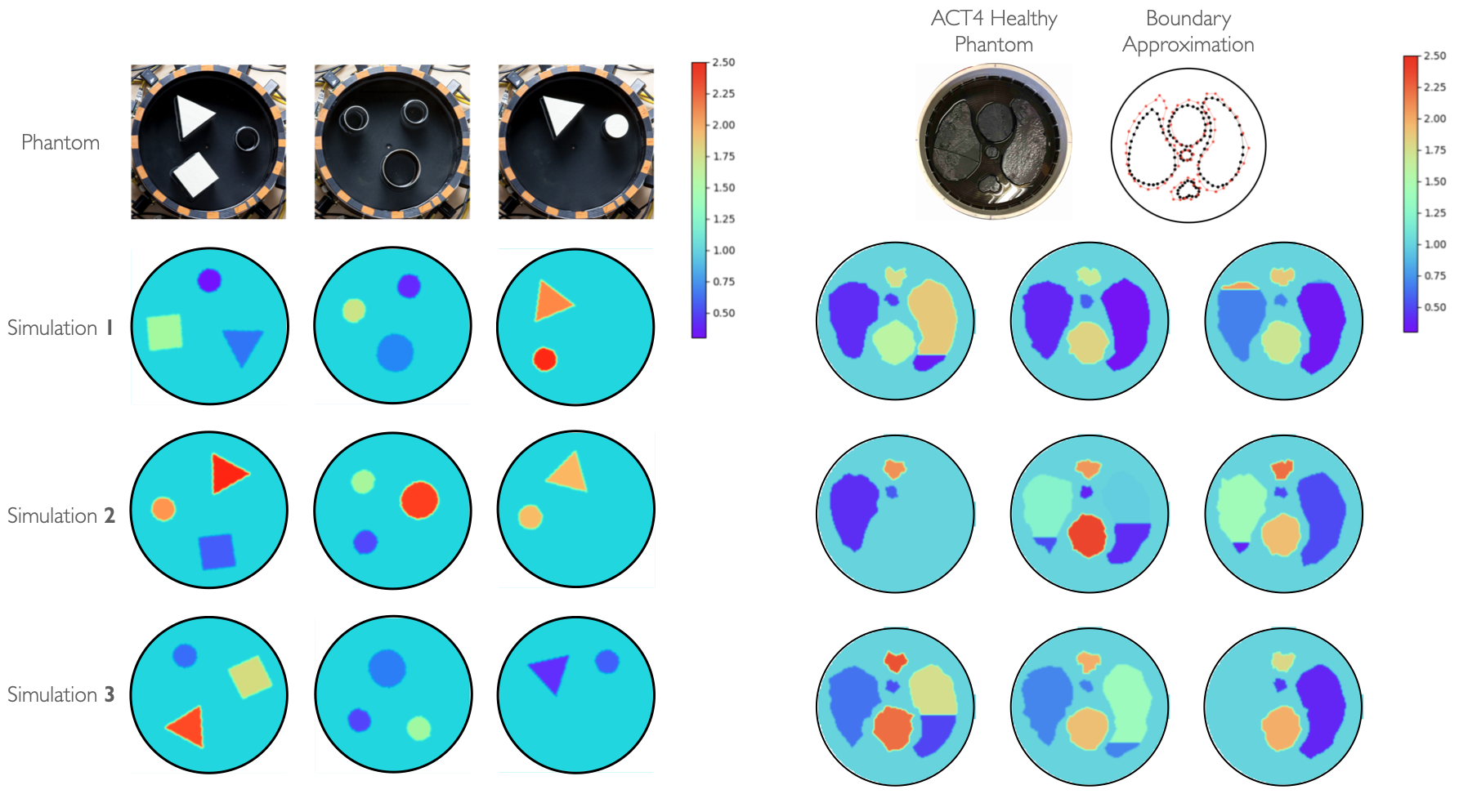} 
   \caption[]{Illustrative examples in KIT4 and ACT4 datasets}
  \label{Dataset}
\end{figure}

Meanwhile, the ACT4 dataset is created similarly to \cite{hamilton2018deep}. Using the ‘HLSA Healthy’ image shown in Fig.~\ref{Dataset} left top, we manually extract the approximated organ boundaries and add a random Gaussian perturbation to the boundary points' locations. Such simulated boundary shapes will not vary too greatly. Next, some organ injuries are simulated by dividing both lungs into two portions with randomly selected horizontal lines. A given probability ($50\%$ chance) determines whether such injuries occur. Besides, we randomly generate conductivity values in the above-divided portions from the uniform distribution $[0.3, 2.5]$. Here, we have simulated a total of \textbf{3200 training samples and 800 testing samples} for each noise level $\delta$. More simulations can be generated for better generalization ability, but this is not considered in the scope of this study. See Fig.~\ref{Dataset} for both illustrative samples of the KIT4 and ACT4 datasets.

\subsection{\textbf{Reconstruction Results}}

Fig.~\ref{KIT4_procedure} and Fig.~\ref{ACT4_procedure} illustrate the workflow of the proposed algorithm for the KIT4 and ACT4 datasets. Each column depicts the processed images, starting with the corresponding ground truth $\sigma(z)$ and ending in the ultimate reconstruction $\tilde{\sigma}(z)$ for each sample. Notably, the second and third rows showcase the low-contrast D-bar reconstruction $\sigma_{r_1}(z)$ and the high-contrast D-bar prediction $\sigma_{r}(z)$, respectively. Specifically, the frequency-enhancing network effectively maps the low-contrast D-bar images to uniform high-contrast D-bar images, significantly improving the visibility and distinction of interior structure. Building on this, the image-calibrating network operates on uniformly high-contrast images, enabling more accurate predictions for EIT reconstruction. 

\begin{figure}[!htp]
  \centering
  \scalebox{1}{
  \begin{tabular}{cccccccccc}
    %\hline
    %Items & Advantages & Disadvantages
    &\scriptsize{Sample 1}&
    \scriptsize{Sample 2}&\scriptsize{Sample 3}&\scriptsize{Sample 4}&
    \scriptsize{Sample 5}
    \\
     \put(-10,15){\rotatebox{90}{\scriptsize{\makecell{Ground \\ Truth}}}}&
    \includegraphics[width=.12\linewidth,height=.12\linewidth]{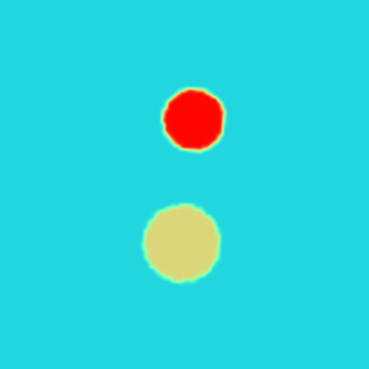}&
    \includegraphics[width=.12\linewidth,height=.12\linewidth]{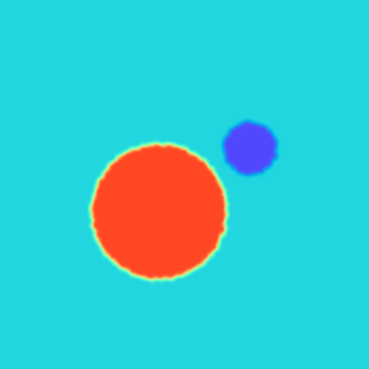}&
    \includegraphics[width=.12\linewidth,height=.12\linewidth]{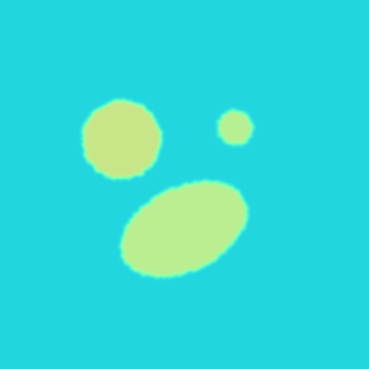}&
    \includegraphics[width=.12\linewidth,height=.12\linewidth]{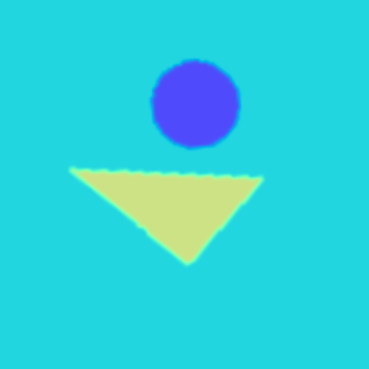}&
    \includegraphics[width=.12\linewidth,height=.12\linewidth]{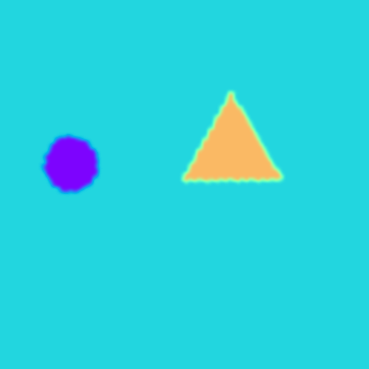}&
    \\ 
    
    \put(-10,5){\rotatebox{90}{\scriptsize{\makecell{Low-Contrast \\ D-bar}}}}&
    \includegraphics[width=.12\linewidth,height=.12\linewidth]{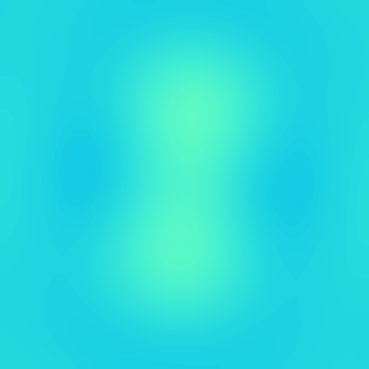}&
    \includegraphics[width=.12\linewidth,height=.12\linewidth]{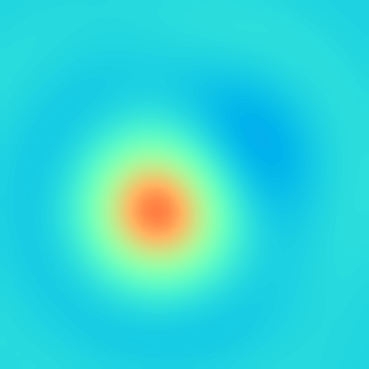}&
    \includegraphics[width=.12\linewidth,height=.12\linewidth]{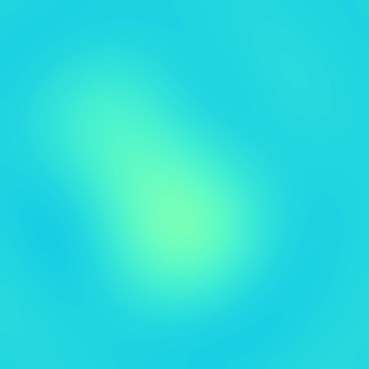}&
    \includegraphics[width=.12\linewidth,height=.12\linewidth]{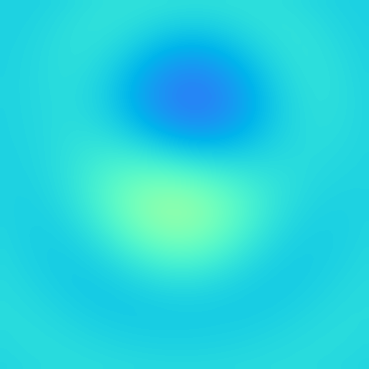}&
    \includegraphics[width=.12\linewidth,height=.12\linewidth]{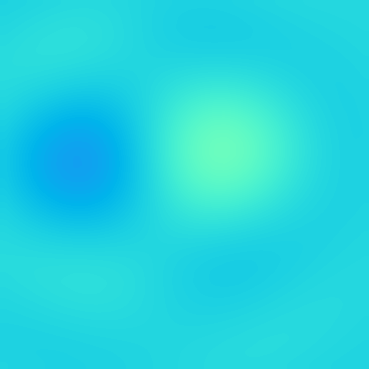}&
    \\ 
    \put(-10,5){\rotatebox{90}{\scriptsize{\makecell{High-Contrast \\ D-bar}}}}&
    \includegraphics[width=.12\linewidth,height=.12\linewidth]{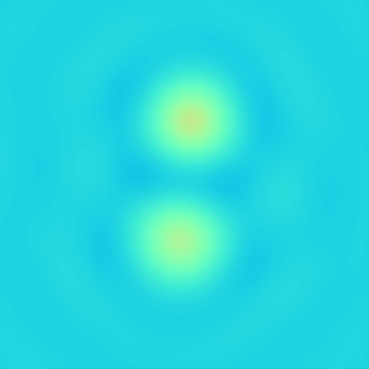}&
    \includegraphics[width=.12\linewidth,height=.12\linewidth]{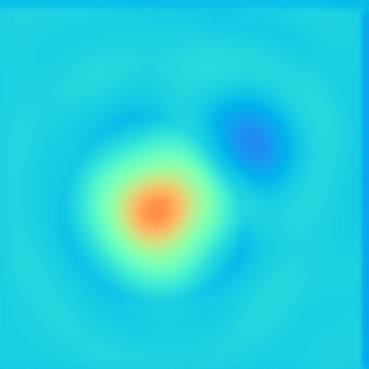}&
    \includegraphics[width=.12\linewidth,height=.12\linewidth]{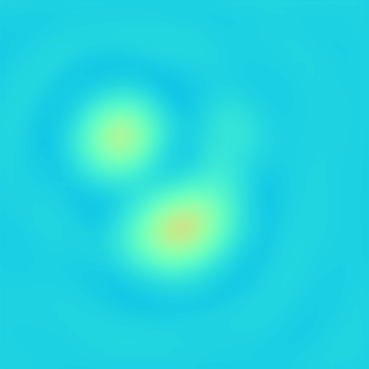}&
    \includegraphics[width=.12\linewidth,height=.12\linewidth]{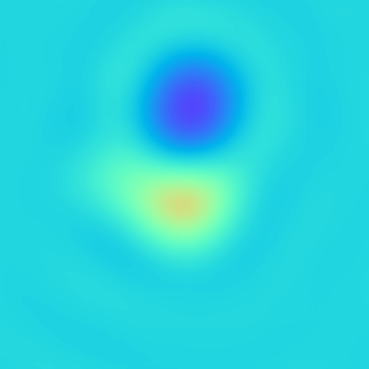}&
    \includegraphics[width=.12\linewidth,height=.12\linewidth]{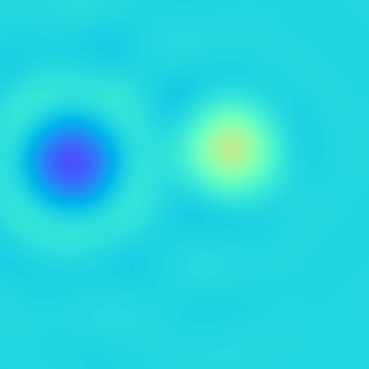}&
    \\ 
    \put(-10,10){\rotatebox{90}{\scriptsize{\makecell{Network \\ Prediction}}}}&
    \includegraphics[width=.12\linewidth,height=.12\linewidth]{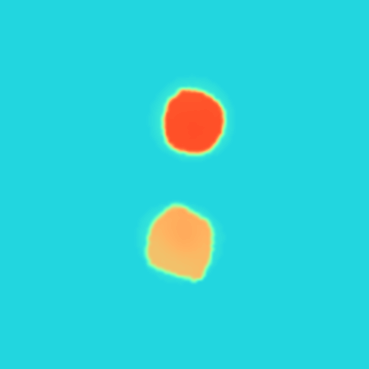}&
    \includegraphics[width=.12\linewidth,height=.12\linewidth]{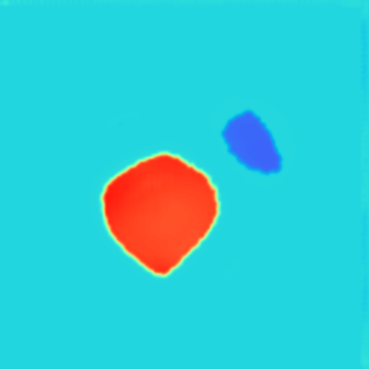}&
    \includegraphics[width=.12\linewidth,height=.12\linewidth]{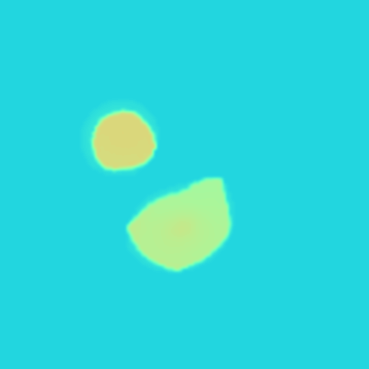}&
    \includegraphics[width=.12\linewidth,height=.12\linewidth]{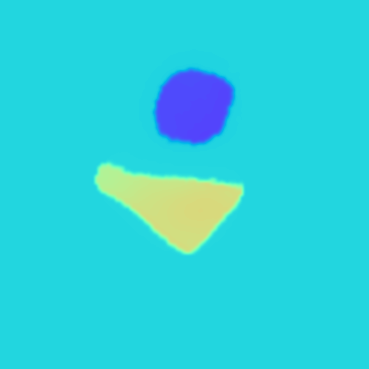}&
    \includegraphics[width=.12\linewidth,height=.12\linewidth]{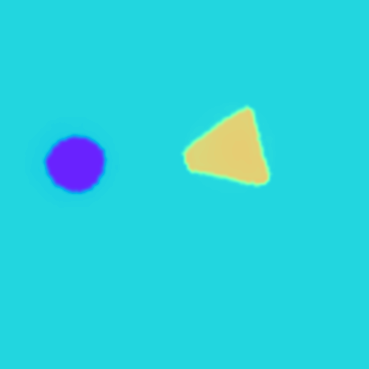}&
  \end{tabular}
  }
  \label{KIT4_procedure}
  
  \centering
		\begin{minipage}[t]{1\linewidth}
	    \hspace{4cm} \includegraphics[width=0.5\linewidth]{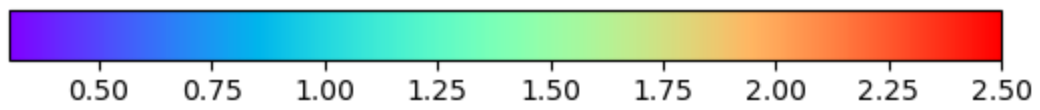}
		\caption{The figure depicts the predicted images during the reconstruction process of the proposed framework for KIT4 dataset. }
		\label{KIT4_procedure}
		\end{minipage}%
	%	\label{KIT4_procedure}
\end{figure}

\begin{figure}[!htp]
  \centering
  \scalebox{1.05}{
  \begin{tabular}{cccccccccc}
    %\hline
    %Items & Advantages & Disadvantages
    &\scriptsize{Sample 1}&\scriptsize{Sample 2}&\scriptsize{Sample 3}&\scriptsize{Sample 4}&\scriptsize{Sample 5}
    \\
    \put(-10,15){\rotatebox{90}{\scriptsize{\makecell{Ground \\ Truth}}}}&
    \includegraphics[width=.12\linewidth,height=.12\linewidth]{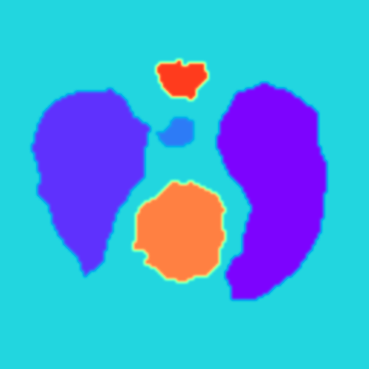}&
    \includegraphics[width=.12\linewidth,height=.12\linewidth]{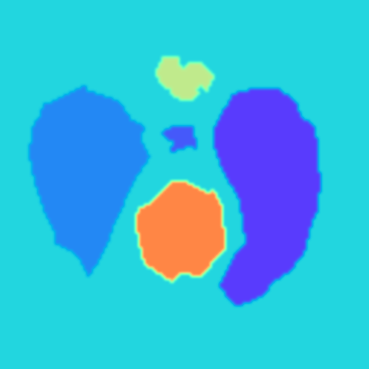}&
    \includegraphics[width=.12\linewidth,height=.12\linewidth]{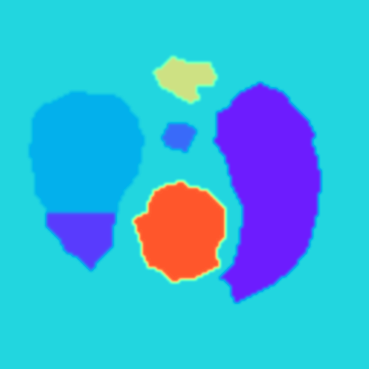}&
    \includegraphics[width=.12\linewidth,height=.12\linewidth]{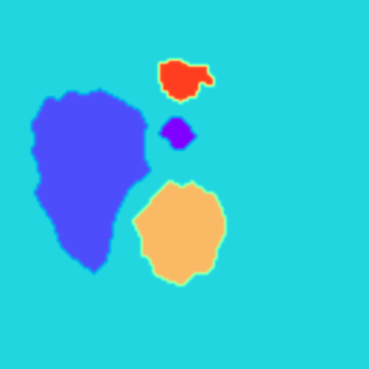}&
    \includegraphics[width=.12\linewidth,height=.12\linewidth]{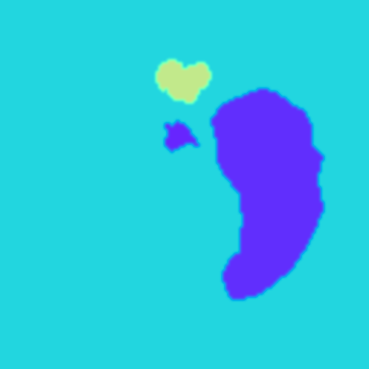}
    \\ 
    \put(-10,5){\rotatebox{90}{\scriptsize{\makecell{Low-Contrast \\ D-bar}}}}&
    \includegraphics[width=.12\linewidth,height=.12\linewidth]{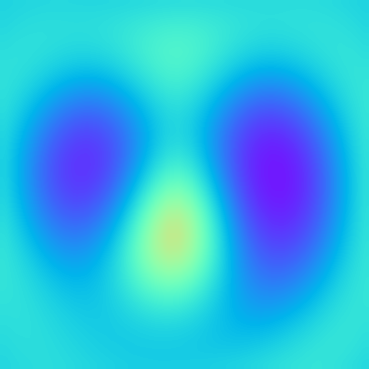}&
    \includegraphics[width=.12\linewidth,height=.12\linewidth]{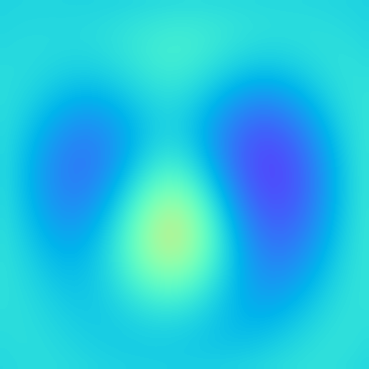}&
    \includegraphics[width=.12\linewidth,height=.12\linewidth]{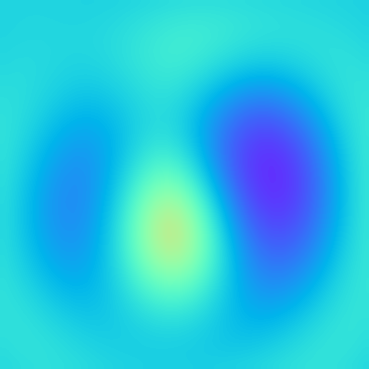}&
    \includegraphics[width=.12\linewidth,height=.12\linewidth]{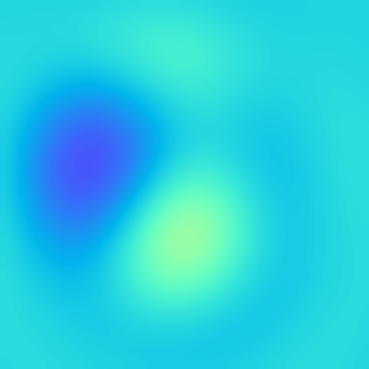}&
    \includegraphics[width=.12\linewidth,height=.12\linewidth]{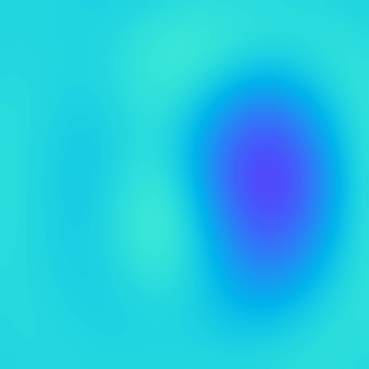}
    \\ 
    \put(-10,5){\rotatebox{90}{\scriptsize{\makecell{High-Contrast \\ D-bar}}}}&
    \includegraphics[width=.12\linewidth,height=.12\linewidth]{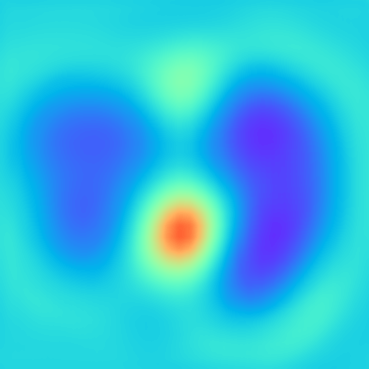}&
    \includegraphics[width=.12\linewidth,height=.12\linewidth]{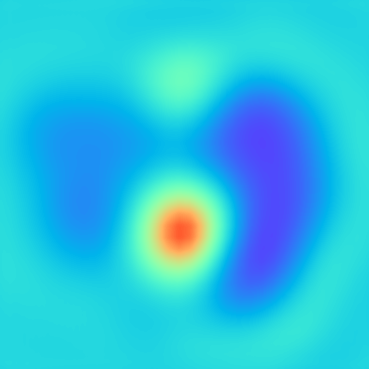}&
    \includegraphics[width=.12\linewidth,height=.12\linewidth]{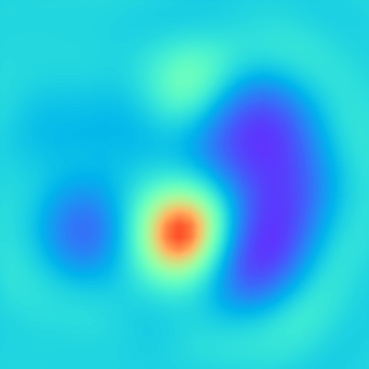}&
    \includegraphics[width=.12\linewidth,height=.12\linewidth]{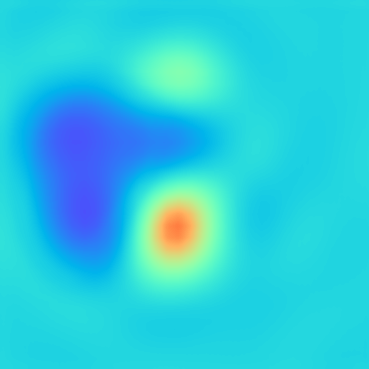}&
    \includegraphics[width=.12\linewidth,height=.12\linewidth]{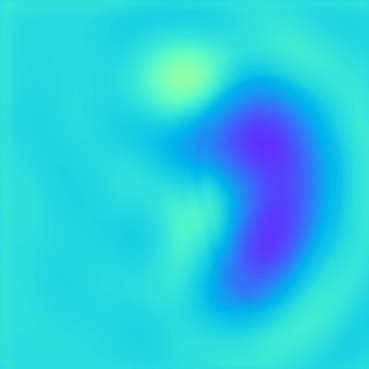}
    \\ 
    
    \put(-10,10){\rotatebox{90}{\scriptsize{\makecell{Network \\ Prediction}}}}&
    \includegraphics[width=.12\linewidth,height=.12\linewidth]{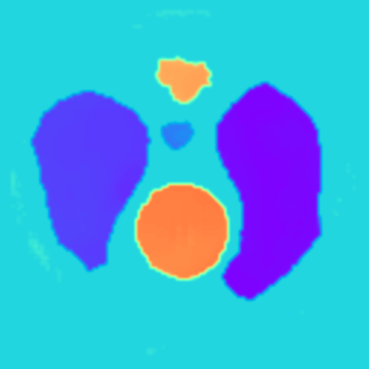}&
    \includegraphics[width=.12\linewidth,height=.12\linewidth]{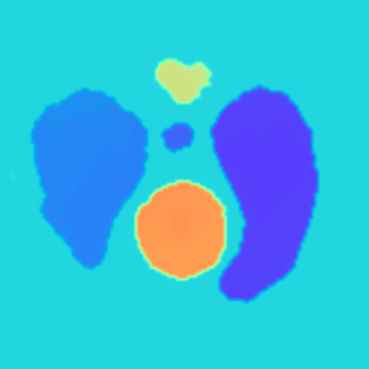}&
    \includegraphics[width=.12\linewidth,height=.12\linewidth]{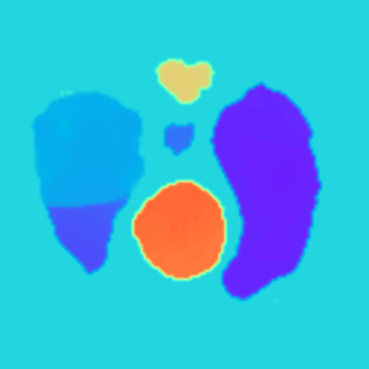}&
    \includegraphics[width=.12\linewidth,height=.12\linewidth]{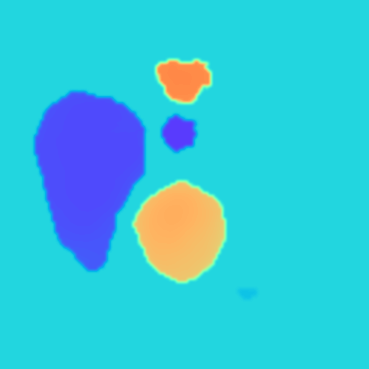}&
    \includegraphics[width=.12\linewidth,height=.12\linewidth]{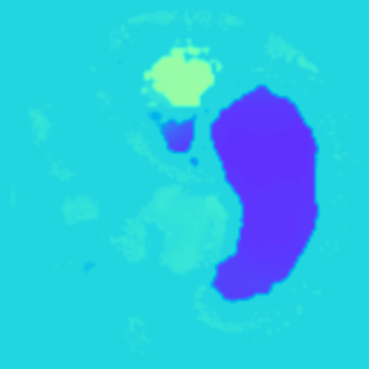}
  \end{tabular}
  }
  
  \label{ACT4_procedure}
  
  \centering
		\begin{minipage}[t]{1\linewidth}
	    \hspace{4.2cm} \includegraphics[width=0.5\linewidth]{Figures/colorbar_s.png}
		\caption{The figure depicts the predicted images during the reconstruction process of the proposed framework for ACT4 dataset. }
		\label{ACT4_procedure}
		\end{minipage}%
\end{figure}

Besides, our proposed method will be compared with the following supervised learning-based methods:  (1) CNN  \cite{hu2019image} (2) Deep D-bar \cite{hamilton2018deep} (3) TSDL \cite{ren2019two} and non-learning methods: (1) GMRES \cite{liu2018pyeit} (2) D-bar \cite{mueller2012linear}. Here, we choose three commonly used metrics including \textbf{PSNR} (Peak Signal-to-Noise Ratio), \textbf{SSIM} (Structural Similarity Indices), and \textbf{RMSE} (Relative Mean Square Error) for quality assessment.

Fig.~\ref{KIT4_comp} and Fig.~\ref{ACT4_comp} show a visual comparison of numerical results for KIT4 and ACT4 datasets, where the proposed method outperforms two non-learning methods in terms of the reconstruction quality: GMRES and D-bar methods. The GMRES algorithm implemented for EIT \cite{liu2018pyeit} is an iterative optimization method to minimize the boundary misfit in PDE-based models, which requires the iterative computation of the associated Jacobian matrices. Besides, this method typically relies on some regularization techniques, such as total variation or Tikhonov regularization, to stabilize the inversion process. The regularized D-bar method \cite{mueller2012linear} applies low-pass truncation on the scattering data for D-bar reconstruction. This approach mitigates the numerical instability and the amplified noise, but also significantly deteriorates reconstruction quality, as high-frequency information for capturing fine details and sharp contrasts is lost.

\begin{figure}[!htp]
  \centering
  \scalebox{0.92}{
  \begin{tabular}{cccccccccc}
    %\hline
    %Items & Advantages & Disadvantages
    &\scriptsize{\makecell{Ground truth}}&\scriptsize{GMRES\cite{liu2018pyeit}}&\scriptsize{Dbar\cite{mueller2012linear}}&\scriptsize{CNN-based \cite{hu2019image}}
    &\scriptsize{TSDL\cite{ren2019two}}&
    \scriptsize{Deep-Dbar\cite{hamilton2018deep}}&\scriptsize{Ours}
    \\ 
    \put(-5,5){\rotatebox{90}{\scriptsize{Sample 1}}}&
    \includegraphics[width=.12\linewidth,height=.12\linewidth]{KIT4/sample_1/gt.png}&
    \includegraphics[width=.12\linewidth,height=.12\linewidth]{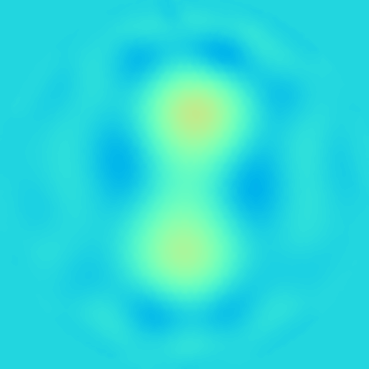}&
    \includegraphics[width=.12\linewidth,height=.12\linewidth]{KIT4/sample_1/dbar.png}&
    \includegraphics[width=.12\linewidth,height=.12\linewidth]{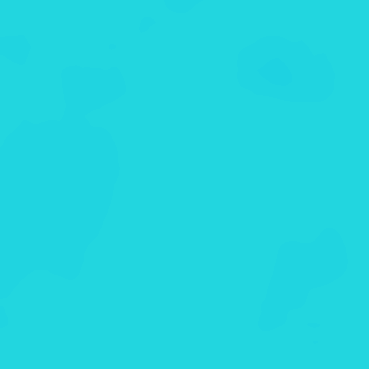}&
    \includegraphics[width=.12\linewidth,height=.12\linewidth]{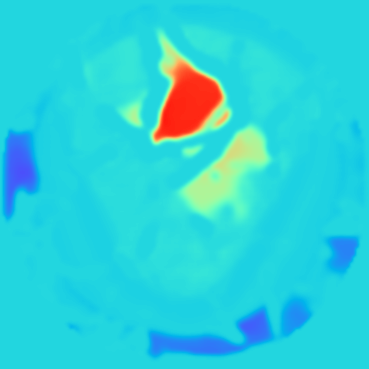}&
    \includegraphics[width=.12\linewidth,height=.12\linewidth]{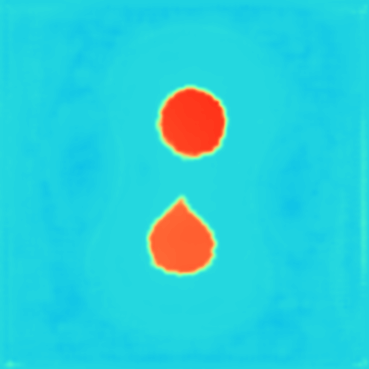}&
    \includegraphics[width=.12\linewidth,height=.12\linewidth]{KIT4/sample_1/ours_image_3.png}
    \\ 
    \put(-5,5){\rotatebox{90}{\scriptsize{Sample 2}}}&
    \includegraphics[width=.12\linewidth,height=.12\linewidth]{KIT4/sample_2/gt.png}&
    \includegraphics[width=.12\linewidth,height=.12\linewidth]{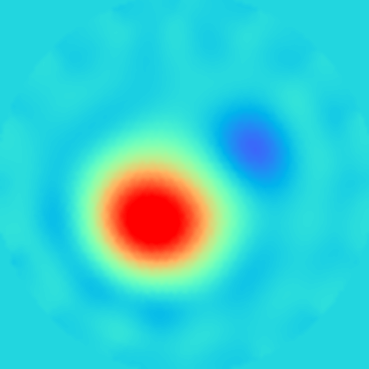}&
    \includegraphics[width=.12\linewidth,height=.12\linewidth]{KIT4/sample_2/dbar.png}&
    \includegraphics[width=.12\linewidth,height=.12\linewidth]{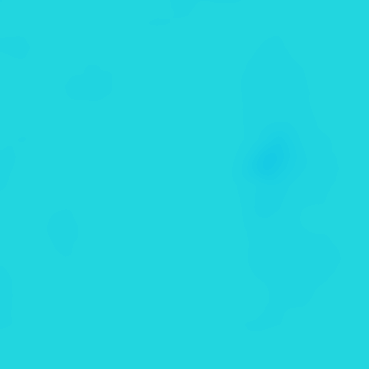}&
    \includegraphics[width=.12\linewidth,height=.12\linewidth]{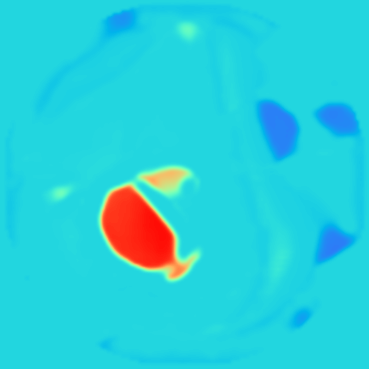}&
    \includegraphics[width=.12\linewidth,height=.12\linewidth]{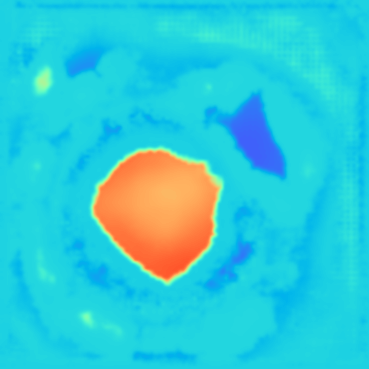}&
    \includegraphics[width=.12\linewidth,height=.12\linewidth]{KIT4/sample_2/ours_image_3.png}
    \\ 
    \put(-5,5){\rotatebox{90}{\scriptsize{Sample 3}}}&
    \includegraphics[width=.12\linewidth,height=.12\linewidth]{KIT4/sample_3/gt.png}&
    \includegraphics[width=.12\linewidth,height=.12\linewidth]{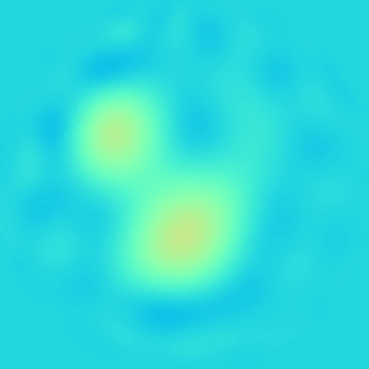}&
    \includegraphics[width=.12\linewidth,height=.12\linewidth]{KIT4/sample_3/dbar.png}&
    \includegraphics[width=.12\linewidth,height=.12\linewidth]{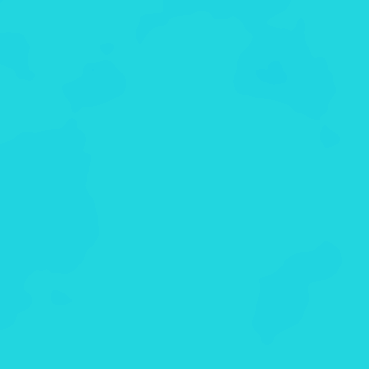}&
    \includegraphics[width=.12\linewidth,height=.12\linewidth]{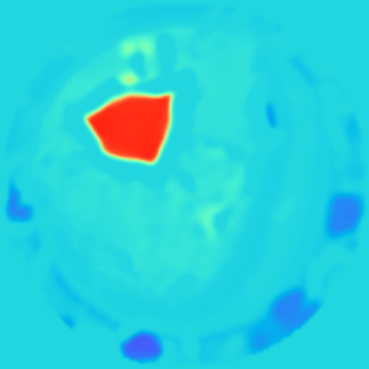}&
    \includegraphics[width=.12\linewidth,height=.12\linewidth]{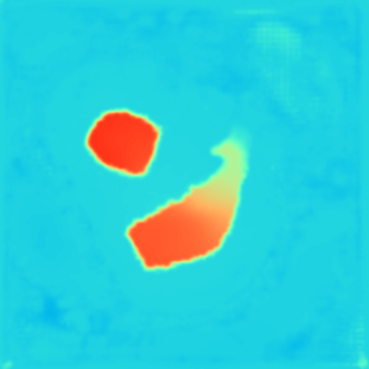}&
    \includegraphics[width=.12\linewidth,height=.12\linewidth]{KIT4/sample_3/ours_image_3.png}
    \\ 
    \put(-5,5){\rotatebox{90}{\scriptsize{Sample 4}}}&
    \includegraphics[width=.12\linewidth,height=.12\linewidth]{KIT4/sample_4/gt.png}&
    \includegraphics[width=.12\linewidth,height=.12\linewidth]{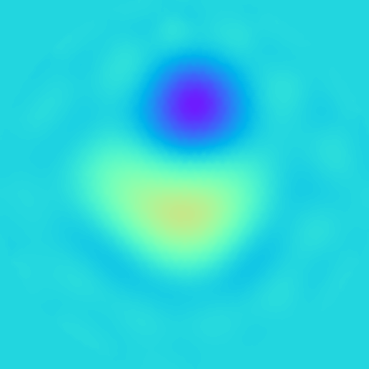}&
    \includegraphics[width=.12\linewidth,height=.12\linewidth]{KIT4/sample_4/dbar.png}&
    \includegraphics[width=.12\linewidth,height=.12\linewidth]{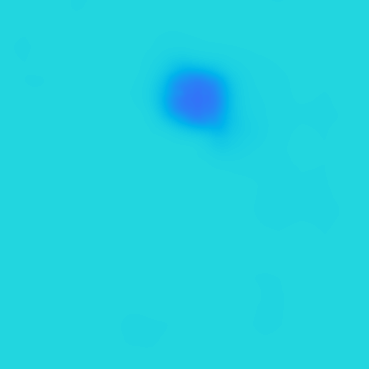}&
    \includegraphics[width=.12\linewidth,height=.12\linewidth]{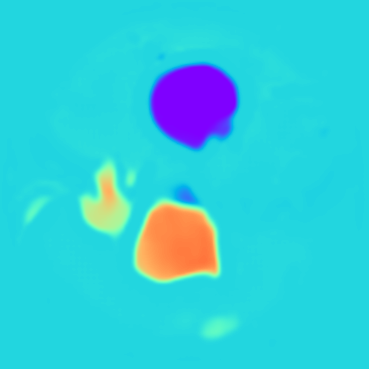}&
    \includegraphics[width=.12\linewidth,height=.12\linewidth]{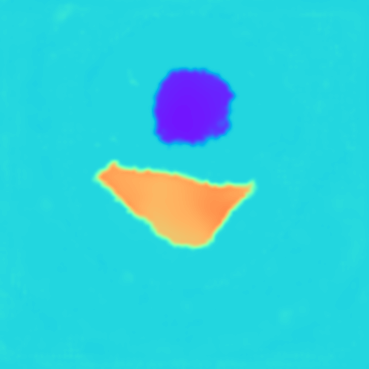}&
    \includegraphics[width=.12\linewidth,height=.12\linewidth]{KIT4/sample_4/ours_image_3.png}
    \\ 
    \put(-5,5){\rotatebox{90}{\scriptsize{Sample 5}}}&
    \includegraphics[width=.12\linewidth,height=.12\linewidth]{KIT4/sample_5/gt.png}&
    \includegraphics[width=.12\linewidth,height=.12\linewidth]{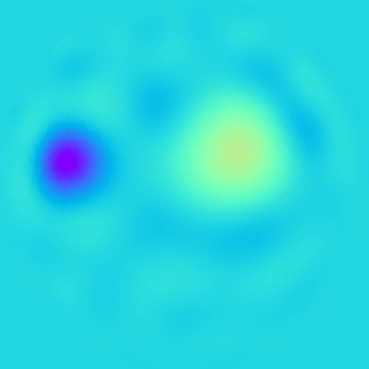}&
    \includegraphics[width=.12\linewidth,height=.12\linewidth]{KIT4/sample_5/dbar.png}&
    \includegraphics[width=.12\linewidth,height=.12\linewidth]{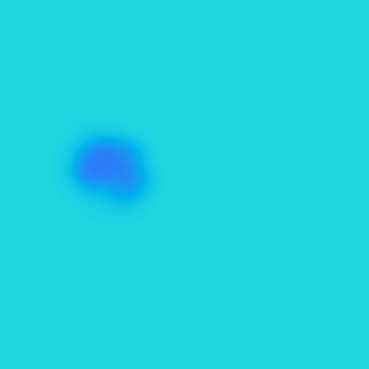}&
    \includegraphics[width=.12\linewidth,height=.12\linewidth]{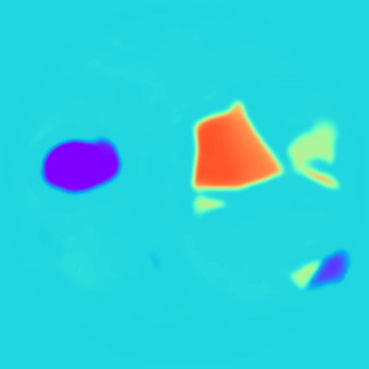}&
    \includegraphics[width=.12\linewidth,height=.12\linewidth]{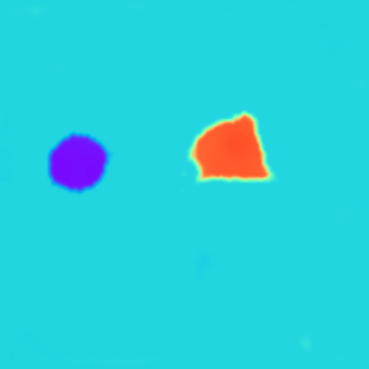}&
    \includegraphics[width=.12\linewidth,height=.12\linewidth]{KIT4/sample_5/ours_image_3.png}
  \end{tabular}
  }
  \label{KIT4_comp}
  \centering
		\begin{minipage}[t]{1\linewidth}
	    \hspace{4.4cm} \includegraphics[width=0.5\linewidth]{Figures/colorbar_s.png}
		\caption{Each row depicts the corresponding ground truth and reconstruction results produced by comparison methods and ours for each KIT4 samples.}
		\label{KIT4_comp}
		\end{minipage}%
\end{figure}

\begin{figure}[!htp]
  \centering
  \scalebox{0.95}{
  \begin{tabular}{cccccccccc}
    %\hline
    %Items & Advantages & Disadvantages
    &\scriptsize{Ground truth}&\scriptsize{GMRES\cite{liu2018pyeit}}&\scriptsize{D-bar\cite{mueller2012linear}}&\scriptsize{CNN-based \cite{hu2019image}}
    &\scriptsize{TSDL\cite{ren2019two}}&
    \scriptsize{Deep D-bar\cite{hamilton2018deep}}&\scriptsize{Ours}
    \\ 
    \put(-5,5){\rotatebox{90}{\scriptsize{Sample 1}}}&
    \includegraphics[width=.12\linewidth,height=.12\linewidth]{ACT4/sample_1/gt.png}&
    \includegraphics[width=.12\linewidth,height=.12\linewidth]{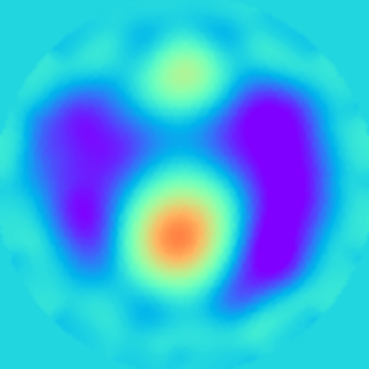}&
    \includegraphics[width=.12\linewidth,height=.12\linewidth]{ACT4/sample_1/dbar.png}&
    \includegraphics[width=.12\linewidth,height=.12\linewidth]{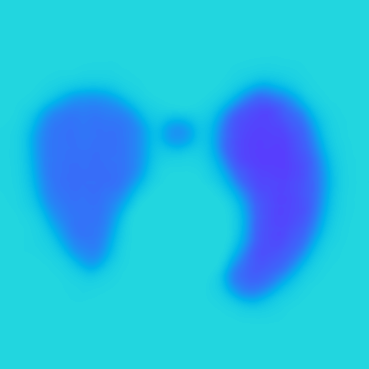}&
    \includegraphics[width=.12\linewidth,height=.12\linewidth]{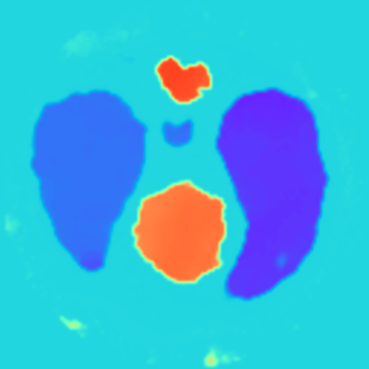}&
    \includegraphics[width=.12\linewidth,height=.12\linewidth]{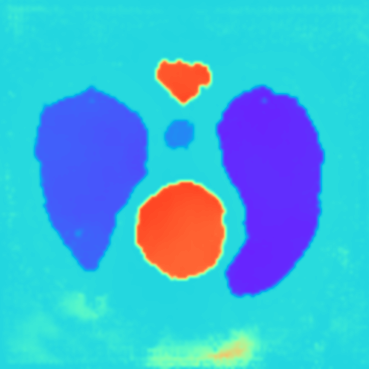}&
    \includegraphics[width=.12\linewidth,height=.12\linewidth]{ACT4/sample_1/ours_image_3.png}

    \\ 
    \put(-5,5){\rotatebox{90}{\scriptsize{Sample 2}}}&
    \includegraphics[width=.12\linewidth,height=.12\linewidth]{ACT4/sample_2/gt.png}&
    \includegraphics[width=.12\linewidth,height=.12\linewidth]{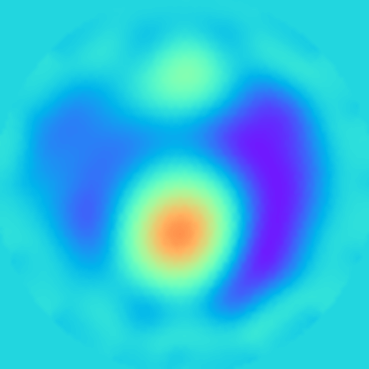}&
    \includegraphics[width=.12\linewidth,height=.12\linewidth]{ACT4/sample_2/dbar.png}&
    \includegraphics[width=.12\linewidth,height=.12\linewidth]{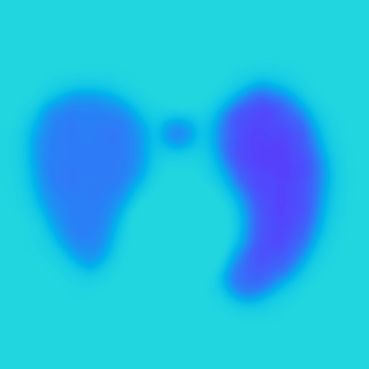}&
    \includegraphics[width=.12\linewidth,height=.12\linewidth]{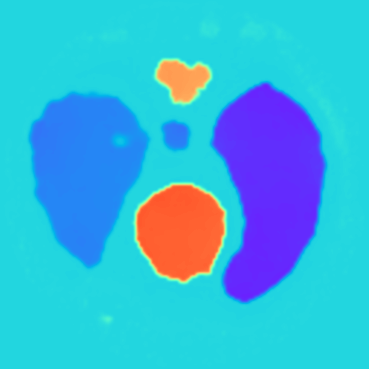}&
    \includegraphics[width=.12\linewidth,height=.12\linewidth]{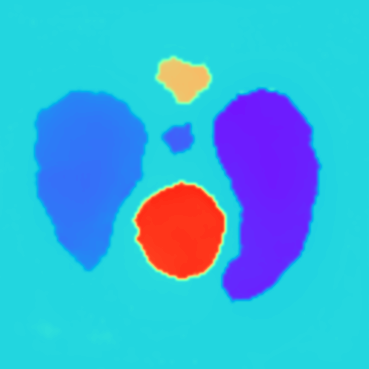}&
    \includegraphics[width=.12\linewidth,height=.12\linewidth]{ACT4/sample_2/ours_image_3.png}
    \\ 
    \put(-5,5){\rotatebox{90}{\scriptsize{Sample 3}}}&
    \includegraphics[width=.12\linewidth,height=.12\linewidth]{ACT4/sample_3/gt.png}&
    \includegraphics[width=.12\linewidth,height=.12\linewidth]{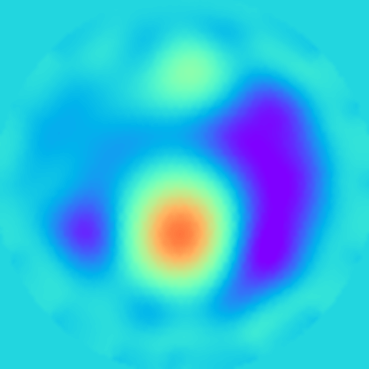}&
    \includegraphics[width=.12\linewidth,height=.12\linewidth]{ACT4/sample_3/dbar.png}&
    \includegraphics[width=.12\linewidth,height=.12\linewidth]{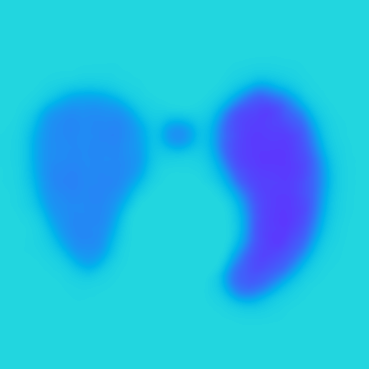}&
    \includegraphics[width=.12\linewidth,height=.12\linewidth]{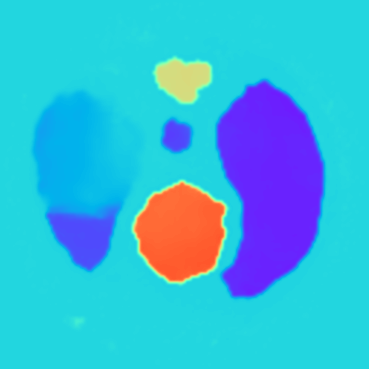}&
    \includegraphics[width=.12\linewidth,height=.12\linewidth]{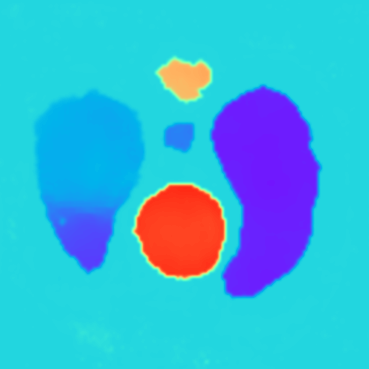}&
    \includegraphics[width=.12\linewidth,height=.12\linewidth]{ACT4/sample_3/ours_image_3.png}
    \\ 
    \put(-5,5){\rotatebox{90}{\scriptsize{Sample 4}}}&
    \includegraphics[width=.12\linewidth,height=.12\linewidth]{ACT4/sample_4/gt.png}&
    \includegraphics[width=.12\linewidth,height=.12\linewidth]{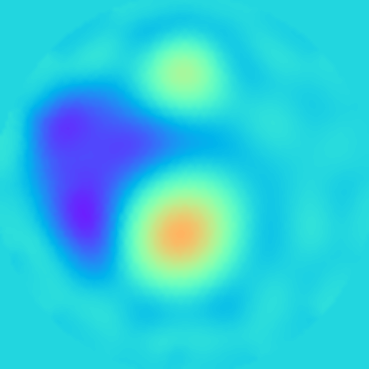}&
    \includegraphics[width=.12\linewidth,height=.12\linewidth]{ACT4/sample_4/dbar.png}&
    \includegraphics[width=.12\linewidth,height=.12\linewidth]{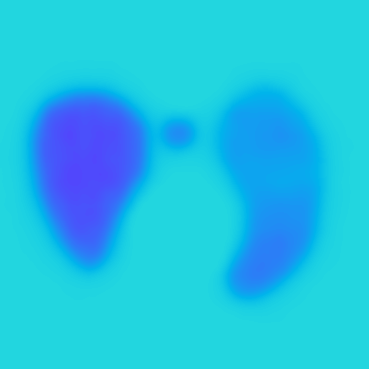}&
    \includegraphics[width=.12\linewidth,height=.12\linewidth]{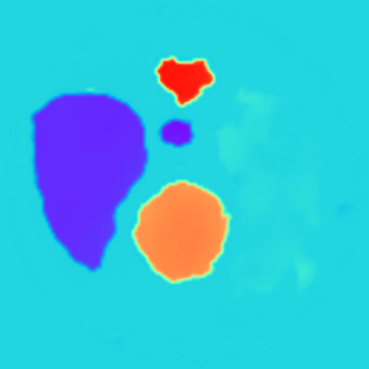}&
    \includegraphics[width=.12\linewidth,height=.12\linewidth]{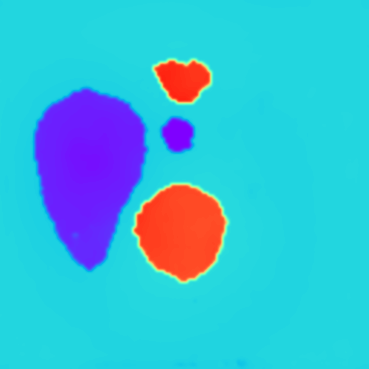}&
    \includegraphics[width=.12\linewidth,height=.12\linewidth]{ACT4/sample_4/ours_image_3.png}
    \\ 
    \put(-5,5){\rotatebox{90}{\scriptsize{Sample 5}}}&
    \includegraphics[width=.12\linewidth,height=.12\linewidth]{ACT4/sample_5/gt.png}&
    \includegraphics[width=.12\linewidth,height=.12\linewidth]{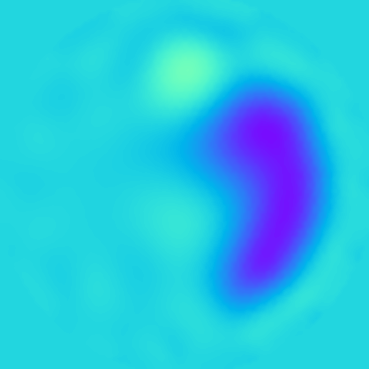}&
    \includegraphics[width=.12\linewidth,height=.12\linewidth]{ACT4/sample_5/dbar.png}&
    \includegraphics[width=.12\linewidth,height=.12\linewidth]{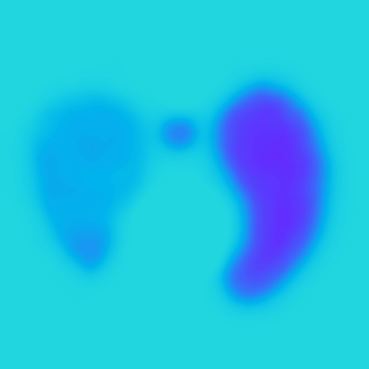}&
    \includegraphics[width=.12\linewidth,height=.12\linewidth]{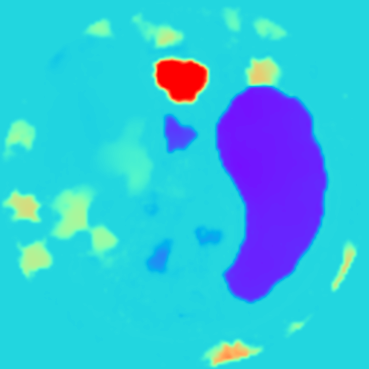}&
    \includegraphics[width=.12\linewidth,height=.12\linewidth]{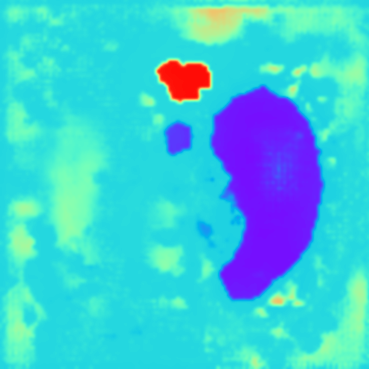}&
    \includegraphics[width=.12\linewidth,height=.12\linewidth]{ACT4/sample_5/ours_image_3.png}
  \end{tabular}
  }
  \label{ACT4_comp}
  \centering
		\begin{minipage}[t]{1\linewidth}
	    \hspace{4.4cm} \includegraphics[width=0.5\linewidth]{Figures/colorbar_s.png}
		\caption{Each row depicts the corresponding ground truth and reconstruction results produced by comparison methods and ours for each ACT4 sample.}
        \label{ACT4_comp}
		\end{minipage}
\end{figure}

Compared to some supervised learning methods, such as the CNN-based approach in \cite{hu2019image} and the TSDL method in \cite{ren2019two}, the proposed method demonstrates superior performance in both visual quality and quantitative metrics.  The CNN-based learning method proposed in \cite{hu2019image} directly approximates the mapping from the measurement operator $\mathbf{L}_{\sigma}^{\delta}$ to the corresponding ground truth conductivity distribution $\sigma(z)$ using a conventional CNN-based architecture for absolute imaging. However, this approach provides only a coarse approximation for negative-value parts as illustrated in the fifth column of Figs.~\ref{KIT4_comp} and \ref{ACT4_comp}. The TSDL method proposed in \cite{ren2019two} features a two-stage architecture comprising pre-reconstruction and post-processing stages. In the first stage, the method back-projects the measurements onto the image space with a pre-computed Jacobian matrix for one-step inversion. In the second stage, a CNN with residual connection will predict the final result from the pre-reconstructed images. The ACT4 experiment, shown in the fourth column of Fig.~\ref{ACT4_comp}, highlights the effectiveness of reconstructing organ phantoms and capturing anatomical structures. However, it produces unexpected artifacts when applied to the KIT4 dataset, as illustrated in Fig.~\ref{KIT4_comp}. In contrast, our proposed method significantly reduces misleading artifacts and preserves the sharp boundaries of inclusions. Furthermore, as indicated in Table.~\ref{quan_result}, it delivers superior results on the KIT4 dataset with comparable results on the ACT4 dataset.

Compared to the single U-Net utilized in the Deep D-bar method \cite{hamilton2018deep}, numerical experiments in Figs.~\ref{KIT4_comp} and \ref{ACT4_comp} demonstrate that our frequency-enhanced learning method achieves more accurate conductivity estimation with less misleading artifacts. For the proposed cascaded architecture, the frequency-enhancement U-Net effectively predicts the frequency-enhanced D-bar image $\sigma_{r_i}(z)$, due to the structural similarity between the input low-contrast D-bar image $\sigma_{r}(z)$ and the output frequency-enhanced D-bar image $\sigma_{r_i}(z)$. Moreover, the image-calibration U-Net serves as a vital post-processing module for prediction. Unlike the Deep D-bar approach post-processing low-contrast D-bar images, the image-calibration U-Net leverages the high-contrast D-bar images as more reliable inputs for prediction. A key advantage of the cascaded structure is its dual-domain utilization of ground truth information, encompassing the non-linear frequency and image domains. It decomposes the overall learning task into hybrid sub-tasks, allowing each network to focus on a specific part during the reconstruction. For the introduced multi-scale structure, we incorporate multi-radii truncated D-bar images into the cascaded learning. This is done for two key purposes: first, the multi-scale structure accommodates multiple selections of the truncation radius compared to the single-scale structure; second, the multi-scale scattering information improves the reconstruction quality, as indicated in Section~\ref{6}. 

\begin{table}
\begin{tabular}{|c|c|c|c|c|c|c|c|c|}
\hline
\multirow{2}{*}{datasets} & \multirow{2}{*}{Metrics} & \multirow{2}{*}{\begin{tabular}[c]{@{}c@{}}Noise\\  Levels\end{tabular}} & \multirow{2}{*}{GMRES} & \multirow{2}{*}{D-bar} & \multirow{2}{*}{CNN} & \multirow{2}{*}{TSDL} & \multirow{2}{*}{\begin{tabular}[c]{@{}c@{}}Deep \\ D-bar\end{tabular}} & \multirow{2}{*}{Ours} \\
                           &                          &                                                                          &                        &                        &                      &                       &                                                                        &                       \\ \hline
\multirow{9}{*}{KIT4}      & \multirow{3}{*}{PSNR$\uparrow$}    & $0.75\%$                                                                 & 24.33                  & 21.04                  & 20.21                & 23.68                 & 24.75                                                                  & \textbf{25.60}               \\ \cline{3-9} 
                           &                          & $0.1\%$                                                                  & 26.70                  & 24.22                  & 22.41                & 24.56                 & 25.78                                                                  & \textbf{27.67}                \\ \cline{3-9} 
                           &                          & $0\%$                                                                    & 27.15                  & 25.65                  & 23.16                & 26.12                 & 26.97                                                                  & \textbf{28.89}            \\ \cline{2-9} 
                           & \multirow{3}{*}{SSIM$\uparrow$}    & $0.75\%$                                                                 & 0.8165                 & 0.7692                 & 0.8025               & 0.8419                & 0.8855                                                                 & \textbf{0.9235}            \\ \cline{3-9} 
                           &                          & $0.1\%$                                                                  & 0.8396                 & 0.7703                 & 0.8065               & 0.8742                & 0.9134                                                                 & \textbf{0.9341}                \\ \cline{3-9} 
                           &                          & $0\%$                                                                    & 0.8714                 & 0.8020                 & 0.8109               & 0.8902                & 0.9176                                                                 & \textbf{0.9455}                \\ \cline{2-9} 
                           & \multirow{3}{*}{RMSE$\downarrow$}    & $0.75\%$                                                                 & 0.1136                & 0.1521                 & 0.2018               & 0.1742                & 0.1213                                                                 & \textbf{0.1096}               \\ \cline{3-9} 
                           &                          & $0.1\%$                                                                  & 0.0868                 & 0.1376                 & 0.1998               & 0.1446                & 0.1021                                                                 & \textbf{0.0812}                \\ \cline{3-9} 
                           &                          & $0\%$                                                                    & 0.0740                 & 0.1143                 & 0.1832               & 0.1321                & 0.0988                                                                 & \textbf{0.0694}                \\ \hline
\multirow{9}{*}{ACT4}      & \multirow{3}{*}{PSNR$\uparrow$}    & $0.75\%$                                                                 & 22.44                  & 19.19                  & 17.33                & 23.24                 & 22.62                                                                  & \textbf{23.63}                 \\ \cline{3-9} 
                           &                          & $0.1\%$                                                                  & 23.12                  & 20.90                  & 18.02                & 24.68                 & 23.57                                                                  & \textbf{25.04}                 \\ \cline{3-9} 
                           &                          & $0\%$                                                                    & 24.46                  & 22.89                  & 18.41               & 25.34                 & 24.43                                                                  & \textbf{26.61}                 \\ \cline{2-9} 
                           & \multirow{3}{*}{SSIM$\uparrow$}    & $0.75\%$                                                                 & 0.6877                 & 0.6556                 & 0.7533               & \textbf{0.8291}                & 0.7949                                                                 & 0.8172                \\ \cline{3-9} 
                           &                          & $0.1\%$                                                                  & 0.7212                 & 0.6754                 & 0.7643               & 0.8457                & 0.8214                                                                 & \textbf{0.8521}                \\ \cline{3-9} 
                           &                          & $0\%$                                                                    & 0.7315                 & 0.6851                 & 0.7706               & \textbf{0.8720}                & 0.8326                                                                 & 0.8657                \\ \cline{2-9} 
                           & \multirow{3}{*}{RMSE$\downarrow$}    & $0.75\%$                                                                 & 0.1588                 & 0.2382                 & 0.2909               & 0.1468                & 0.1578                                                                 & \textbf{0.1445}                \\ \cline{3-9} 
                           &                          & $0.1\%$                                                                  & 0.1321                 & 0.2106                 & 0.2743               & 0.1280                & 0.1499                                                                 & \textbf{0.1223}               \\ \cline{3-9} 
                           &                          & $0\%$                                                                    & 0.1286                 & 0.2047                 & 0.2757               & 0.1193                & 0.1282                                                                 & \textbf{0.1049}                \\ \hline
\end{tabular}
\caption{Quantitative comparison for KIT4 and ACT4 datasets to PSNR, SSIM and RMSE metrics}
\label{quan_result}
\end{table}

\subsection{\textbf{Computational Efficiency}}
Another significant motivation for the proposed method lies in its ability 
for computational efficiency. For a comprehensive comparison, we compare the training/testing time of various methods and GFLOPs/parameters of different networks in Table \ref{timetable}.

\begin{table}[htbp]
\centering
\scalebox{0.92}{
\begin{tabular}{|c|c|c|c|c|c|}
\hline
\multicolumn{1}{|l|}{} & GMRES \cite{liu2018pyeit} & D-bar \cite{mueller2012linear}& TSDL \cite{ren2019two} & Deep D-bar \cite{hamilton2018deep} & Ours \\ \hline
Training Time       & ---     & ---     & 40min    & 40min   & 1h    \\ \hline
Testing Time       & 305.02s    &  54.86s     & 0.025s   &  0.024s + 0.018s &  0.045s + 0.018s\\ \hline
GFLOPs(G)  & ---     & ---   &  5.03   &  1.63   &  3.43  \\ \hline
Params(M)  & ---     & ---   &  0.30   &  1.14   &  2.56  \\ \hline
\end{tabular}
}
\caption{Comparison of training/testing time of various methods and GFLOPs/parameters of different networks.}
\label{timetable}
\end{table}

%Notably, evaluating the proposed cascaded U-Nets is highly computationally efficient on GPUs. As discussed in Subsection \ref{Num_imp_D_bar}, the D-bar reconstruction process can also be substantially accelerated on GPUs. Consequently, our method holds strong potential for real-time EIT reconstruction.

The GMRES algorithm in \cite{liu2018pyeit} requires recursively solving non-linear PDE systems and linear adjoint equations under different boundary conditions. As a result, its computational demand is typically higher than other direct methods, primarily due to the intensive computations needed for the associated Jacobian matrices. The D-bar method in \cite{mueller2012linear} involves hundreds or thousands of calculations of the numerical CGO solutions and the GMRES iteration algorithm for solving D-bar equations. Here, the original D-bar method uses the CPU-based GMRES iterative algorithm for solving D-bar equations. In contrast, the Deep D-bar and our proposed method will use the GPU-based fixed-point iterative algorithm to obtain real-time D-bar reconstructions. The TSDL method in \cite{ren2019two} follows the difference-imaging strategy with a pre-computed Jacobian matrix to project the measurements back onto the image space. Therefore, the training and testing time for TSDL is relatively fast because its reconstruction process is network-based, without time-consuming forward PDE computations or adjoint equation calculations. 
Here, we account for the testing time required for the U-Net post-processing and the accelerated D-bar method. In Table \ref{timetable}, for example, $0.024s$ and $0.018s$ refer to the average times for U-Net post-processing and the fixed-point iteration on the GPU, respectively. In our learning approach, the training process depends on the real-time computation of frequency-enhanced D-bar images for different truncation radii. This step is significantly accelerated by the fixed-point iteration on the GPU, resulting in no substantial increase in training time compared to the Deep D-bar method.

Here, all non-linear PDE solvers are executed on the Intel-Xeon-Platinum-8358P CPU, while the deep learning part and the GPU-based D-bar solver are performed on the NVIDIA-A100 GPU.

\subsection{Out-of-Distribution Testing}

We utilize the models trained on the simulated data from the continuum model to evaluate their performance for more realistic testing data.  Since the KIT4 and ACT4 systems employ electrodes attached along the boundary for data acquisition, we consider the measurement data based on the complete electrode model rather than the continuum model, using $32$ electrodes as described in Section~\ref{4.1}.

\begin{figure}[!htp]
  \centering
  \scalebox{1}{
  \begin{tabular}{ccccccc}
    %\hline
    %Items & Advantages & Disadvantages
    &\scriptsize{Experiment}&\scriptsize{Ground Truth}&\scriptsize{D-bar\cite{mueller2012linear}}&
    \scriptsize{Deep D-bar\cite{hamilton2018deep}}&\scriptsize{Ours}
    \\
    \put(-0,5){\rotatebox{90}{\scriptsize{Phantom 1.1}}}&
    \includegraphics[width=.12\linewidth,height=.12\linewidth]{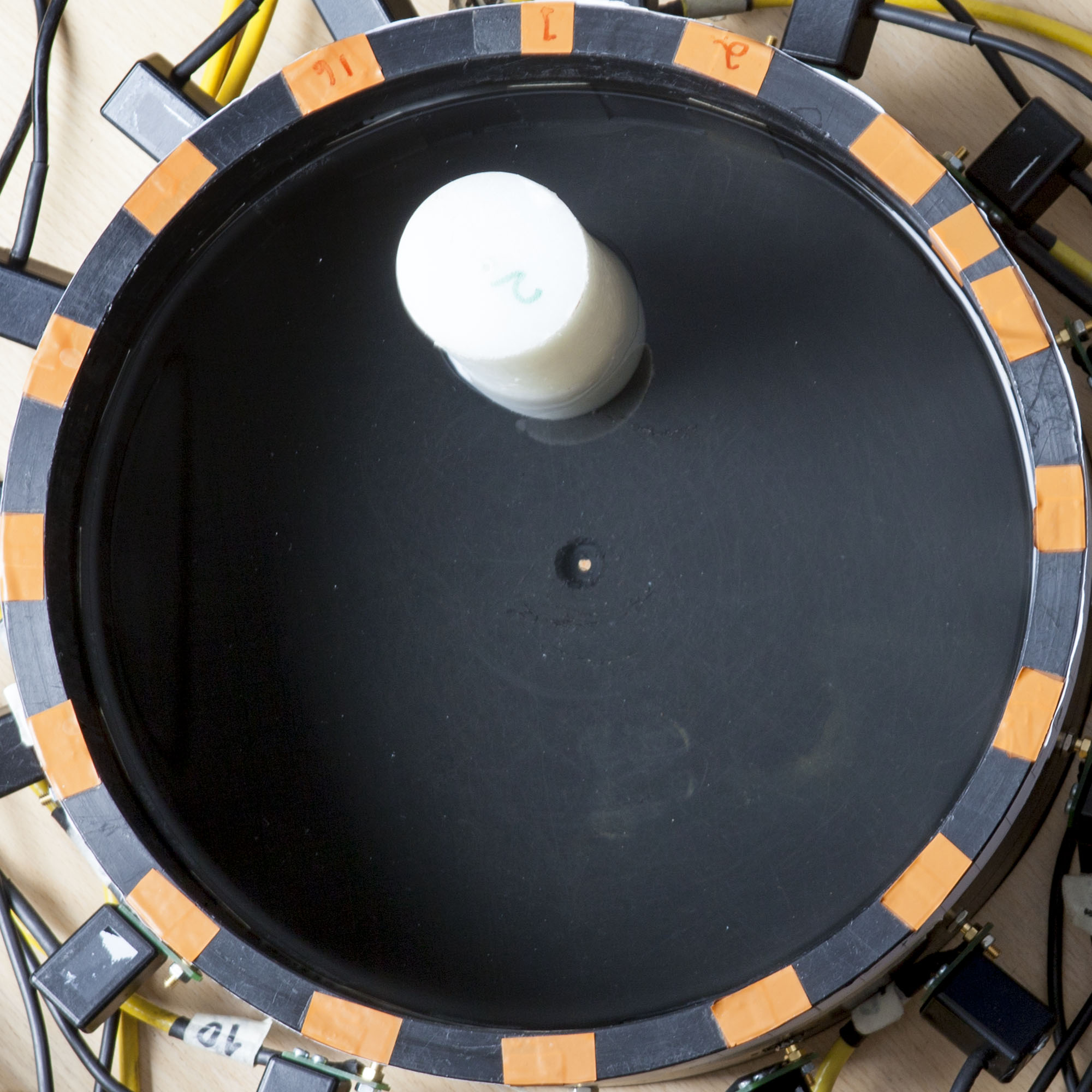}&
    \includegraphics[width=.12\linewidth,height=.12\linewidth]{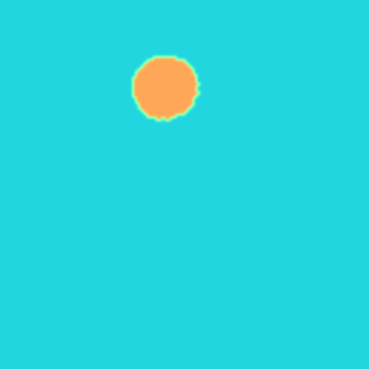}&
    \includegraphics[width=.12\linewidth,height=.12\linewidth]{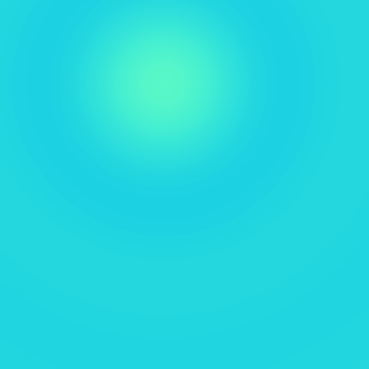}&
    \includegraphics[width=.12\linewidth,height=.12\linewidth]{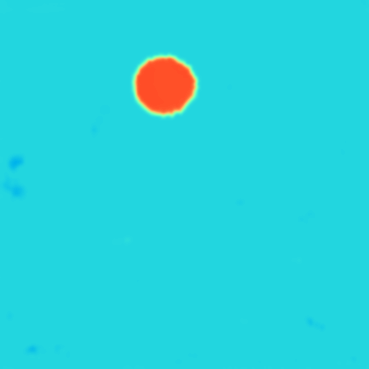}&
    \includegraphics[width=.12\linewidth,height=.12\linewidth]{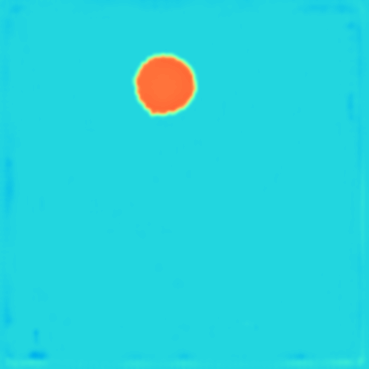}
    \\ 
    \put(-0,5){\rotatebox{90}{\scriptsize{Phantom 2.2}}}&
    \includegraphics[width=.12\linewidth,height=.12\linewidth]{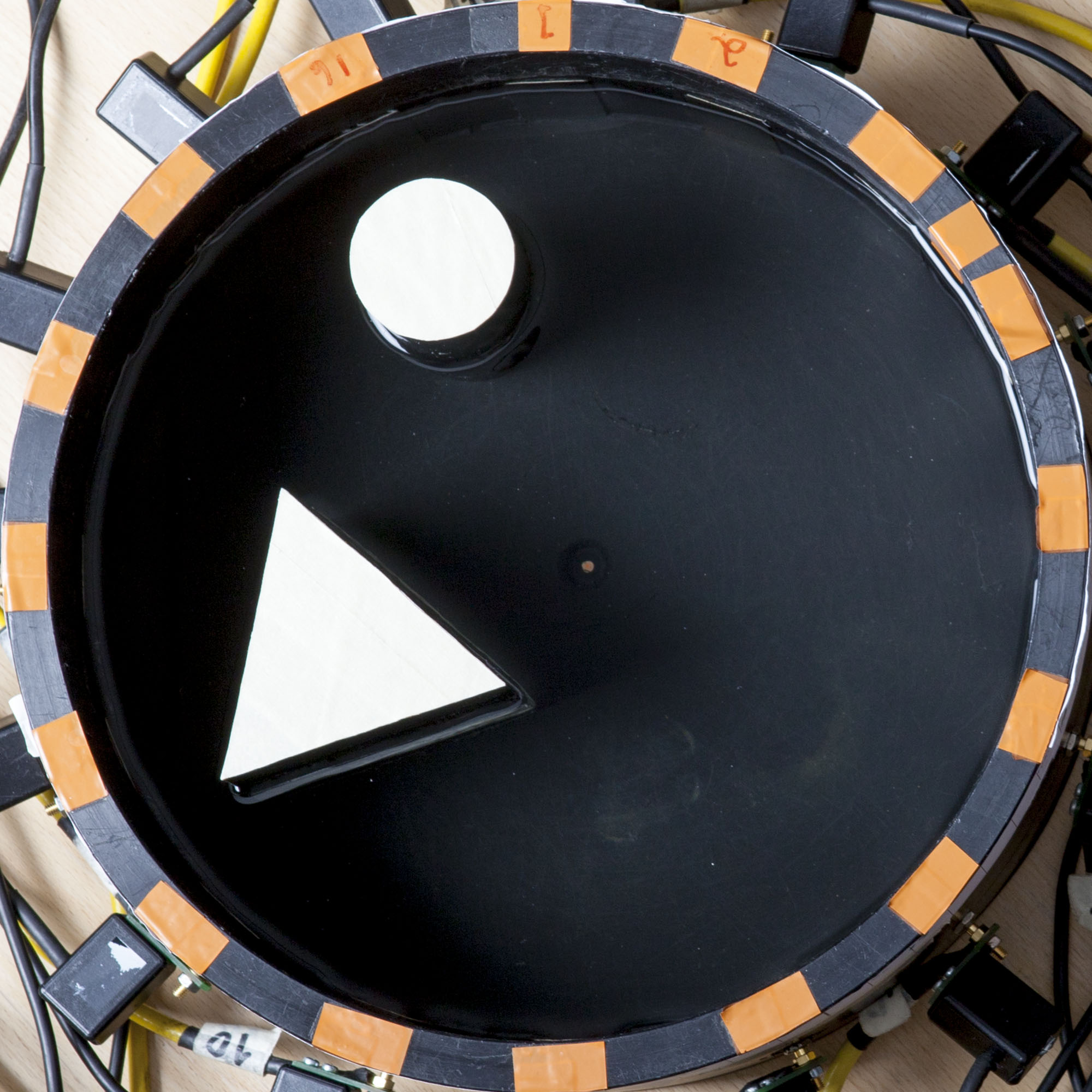}&
    \includegraphics[width=.12\linewidth,height=.12\linewidth]{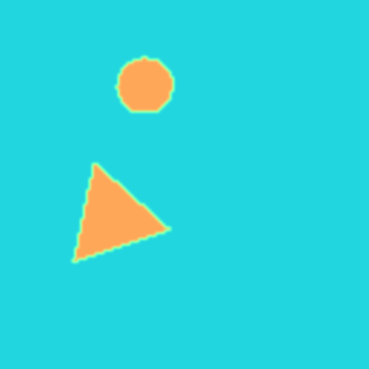}&
    \includegraphics[width=.12\linewidth,height=.12\linewidth]{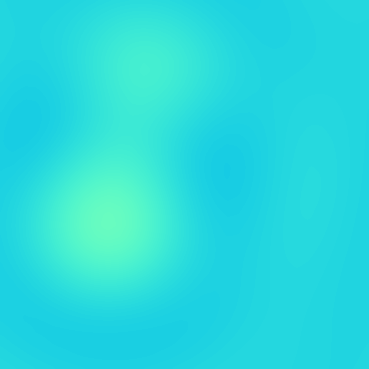}&
    \includegraphics[width=.12\linewidth,height=.12\linewidth]{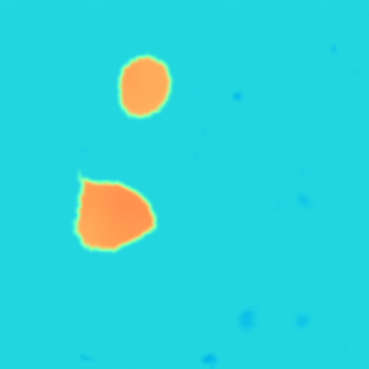}&
    \includegraphics[width=.12\linewidth,height=.12\linewidth]{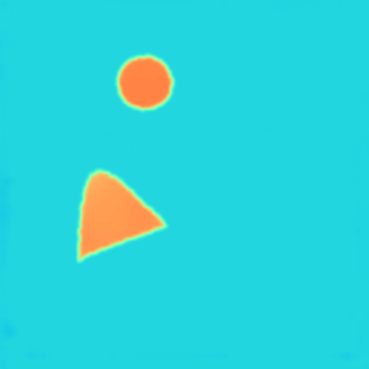}
    \\ 
    \put(-0,5){\rotatebox{90}{\scriptsize{Phantom 3.4}}}&
    \includegraphics[width=.12\linewidth,height=.12\linewidth]{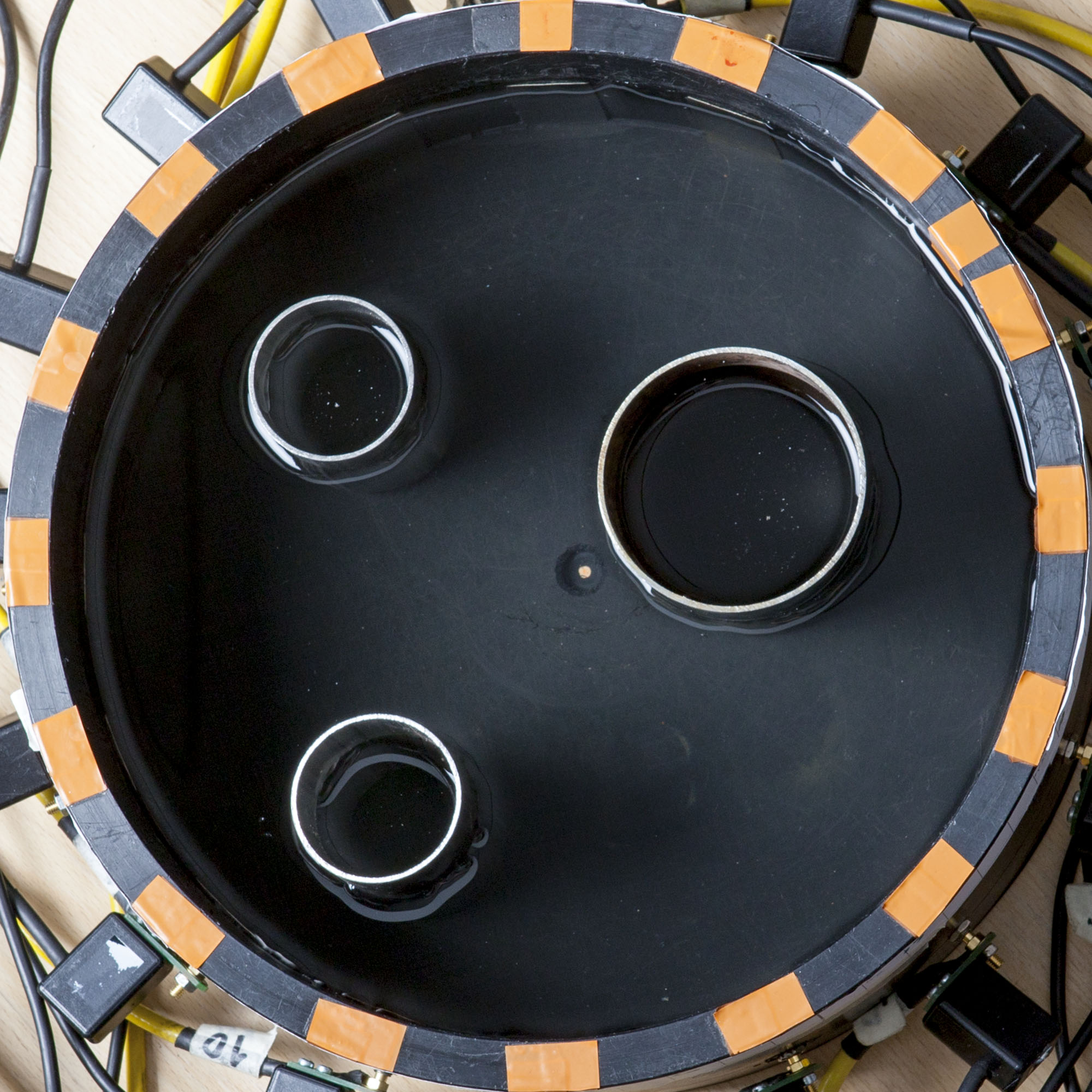}&
    \includegraphics[width=.12\linewidth,height=.12\linewidth]{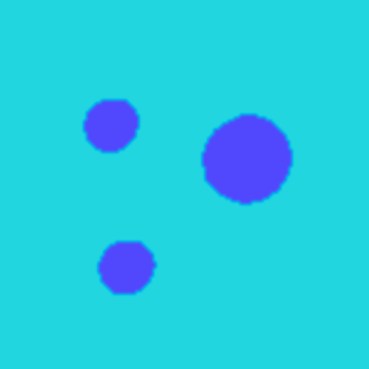}&
    \includegraphics[width=.12\linewidth,height=.12\linewidth]{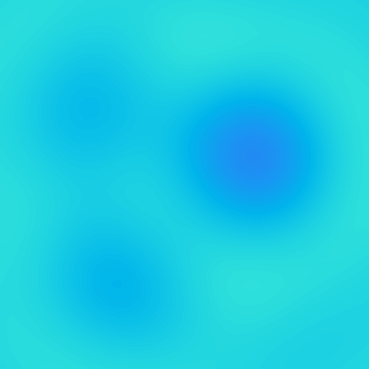}&
    \includegraphics[width=.12\linewidth,height=.12\linewidth]{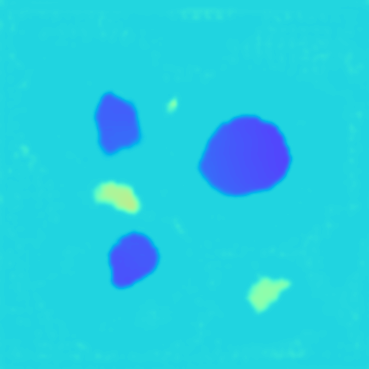}&
    \includegraphics[width=.12\linewidth,height=.12\linewidth]{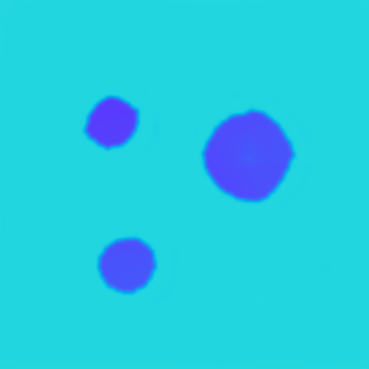}
    \\ 
    \put(-0,5){\rotatebox{90}{\scriptsize{Phantom 4.4}}}&
    \includegraphics[width=.12\linewidth,height=.12\linewidth]{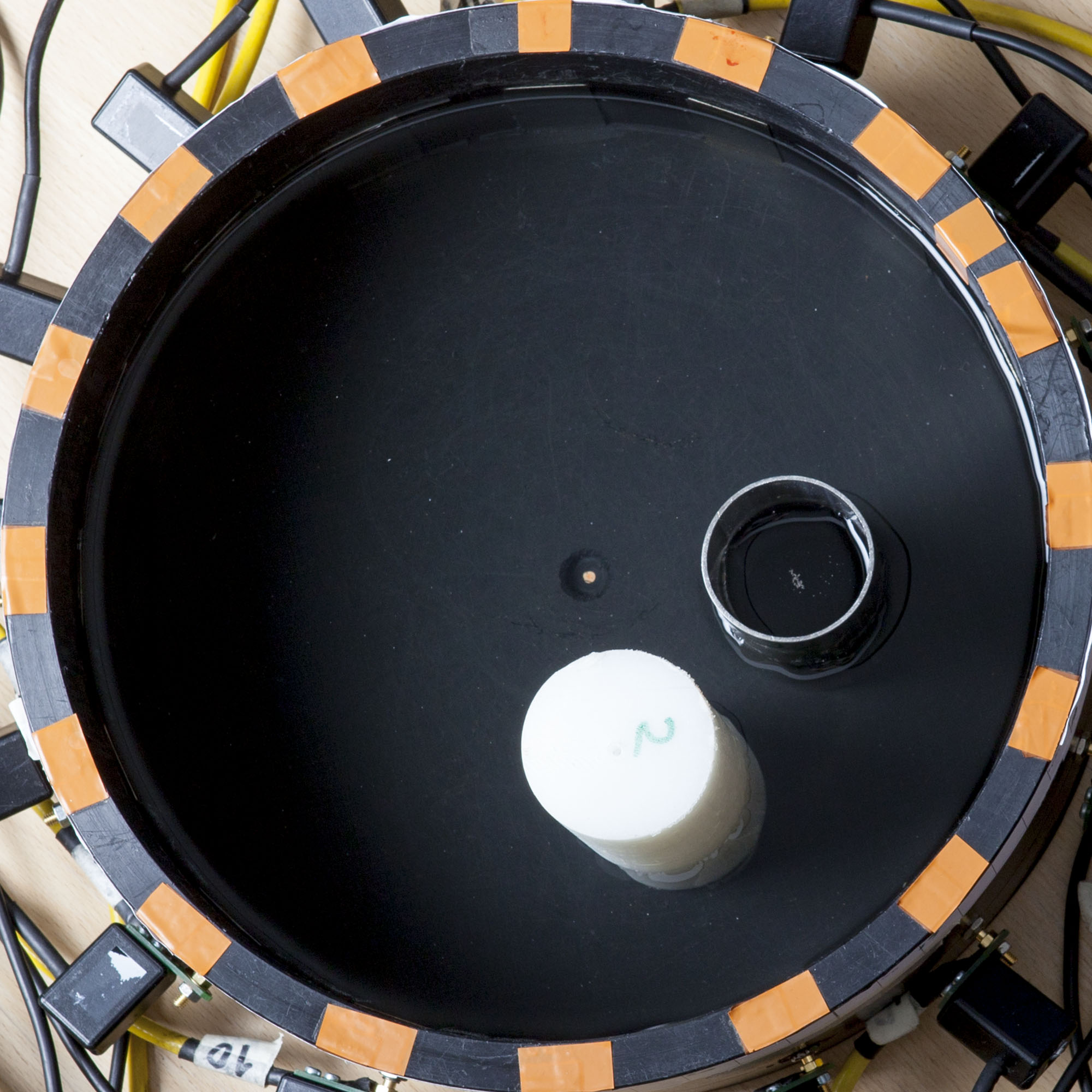}&
    \includegraphics[width=.12\linewidth,height=.12\linewidth]{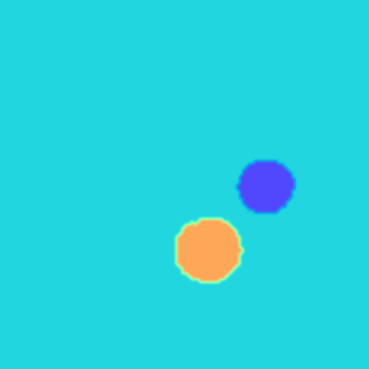}&
    \includegraphics[width=.12\linewidth,height=.12\linewidth]{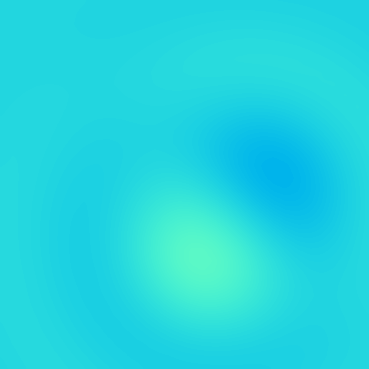}&
    \includegraphics[width=.12\linewidth,height=.12\linewidth]{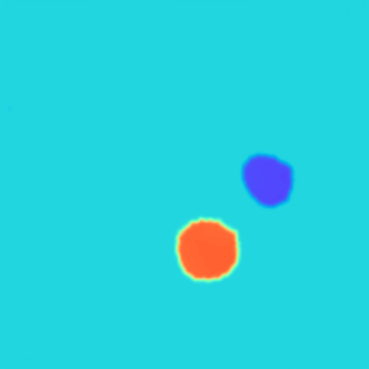}&
    \includegraphics[width=.12\linewidth,height=.12\linewidth]{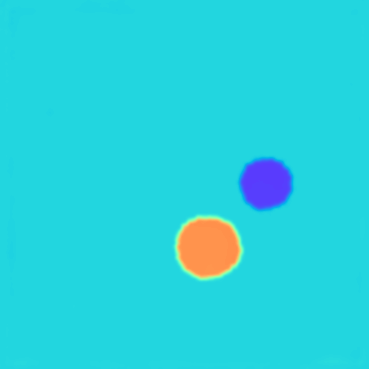}
  \end{tabular}
  } 
  \centering
		\begin{minipage}[t]{1\linewidth}
	    \hspace{4.4cm} \includegraphics[width=0.5\linewidth]{Figures/colorbar_s.png}
		\caption{In the KIT4 system, the measurements are collected from the complete electrode model. The results are produced from models trained on the simulated data from the continuum model.}
        \label{real_world_KIT4}
		\end{minipage}
\end{figure}

\begin{figure}[!htp]
  \centering
  \scalebox{1}{
  \begin{tabular}{cccccc}
    %\hline
    %Items & Advantages & Disadvantages
    &\scriptsize{Ground truth}&\scriptsize{Ground Truth}&\scriptsize{D-bar\cite{mueller2012linear}}&
    \scriptsize{Deep D-bar\cite{hamilton2018deep}}&\scriptsize{Ours}
    \\
    \put(-0,10){\rotatebox{90}{\scriptsize{Healthy}}}&
    \includegraphics[width=.12\linewidth,height=.12\linewidth]{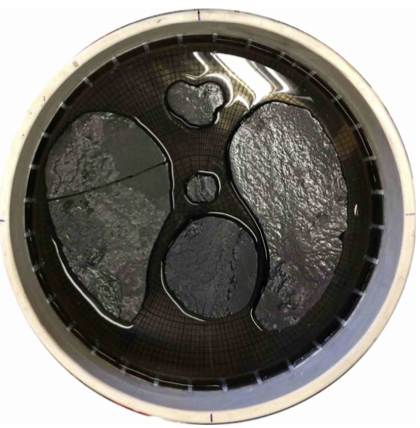}&
    \includegraphics[width=.12\linewidth,height=.12\linewidth]{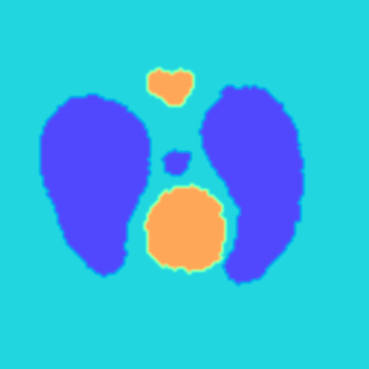}&
    \includegraphics[width=.12\linewidth,height=.12\linewidth]{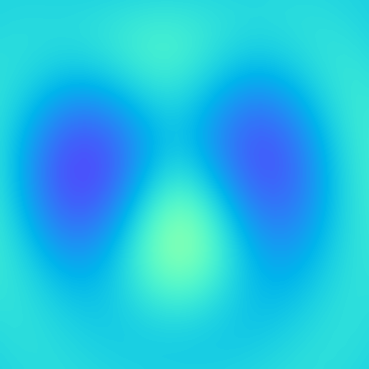}&
    \includegraphics[width=.12\linewidth,height=.12\linewidth]{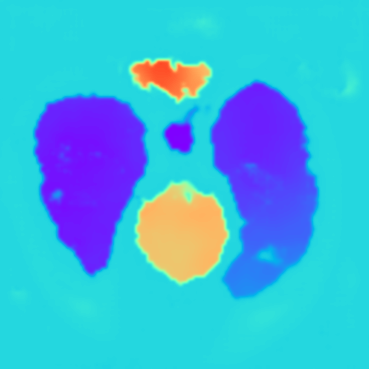}&
    \includegraphics[width=.12\linewidth,height=.12\linewidth]{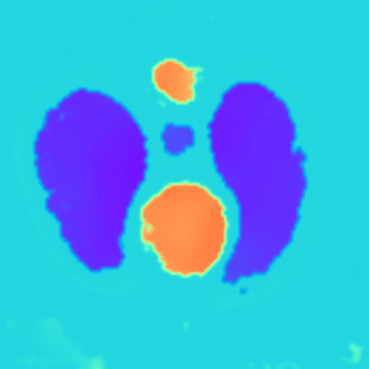}
    \\ 
    \put(-0,10){\rotatebox{90}{\scriptsize{Injury 1}}}&
    \includegraphics[width=.12\linewidth,height=.12\linewidth]{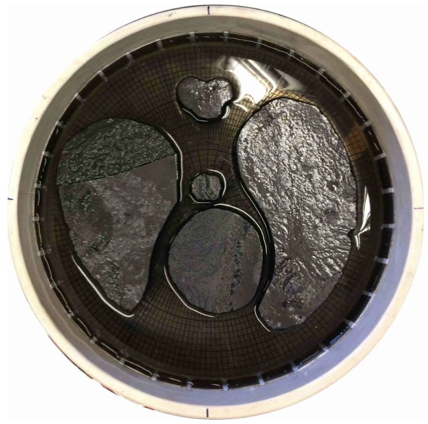}&
    \includegraphics[width=.12\linewidth,height=.12\linewidth]{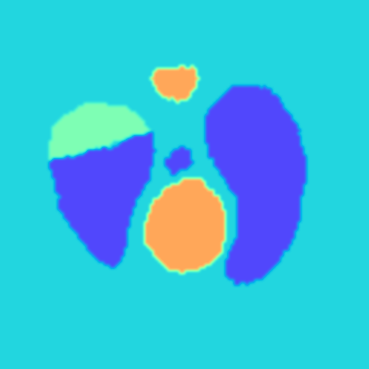}&
    \includegraphics[width=.12\linewidth,height=.12\linewidth]{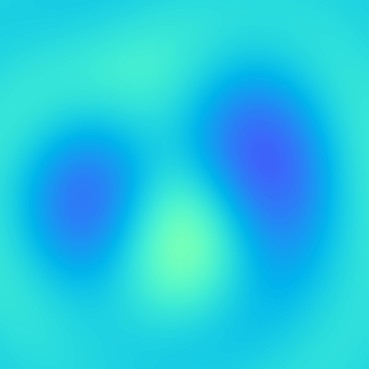}&
    \includegraphics[width=.12\linewidth,height=.12\linewidth]{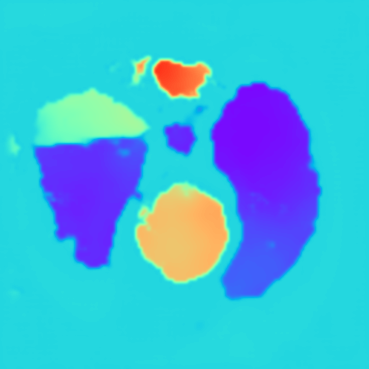}&
    \includegraphics[width=.12\linewidth,height=.12\linewidth]{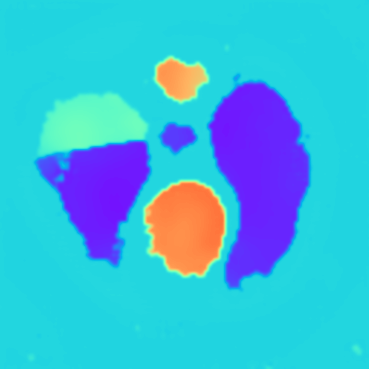}
    \\ 
    \put(-0,10){\rotatebox{90}{\scriptsize{Injury 2}}}&
    \includegraphics[width=.12\linewidth,height=.12\linewidth]{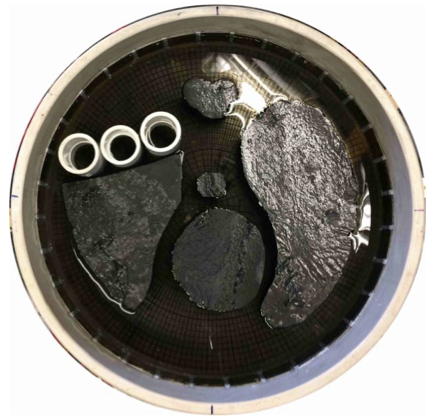}&
    \includegraphics[width=.12\linewidth,height=.12\linewidth]{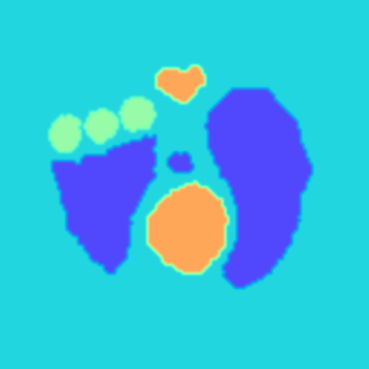}&
    \includegraphics[width=.12\linewidth,height=.12\linewidth]{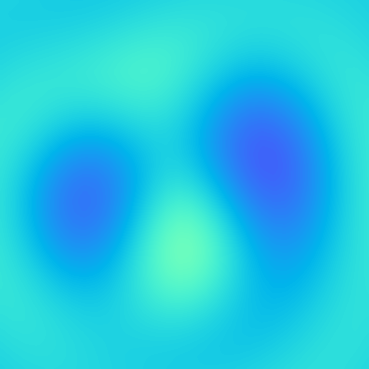}&
    \includegraphics[width=.12\linewidth,height=.12\linewidth]{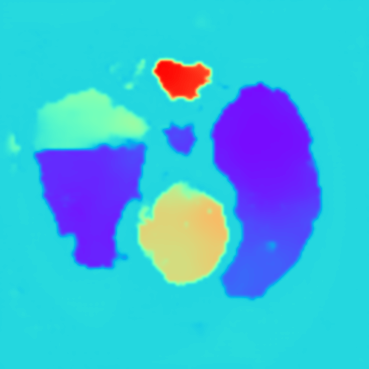}&
    \includegraphics[width=.12\linewidth,height=.12\linewidth]{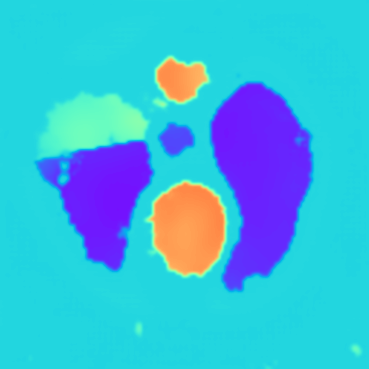}
    \\ 
    \put(-0,10){\rotatebox{90}{\scriptsize{Injury 3}}}&
    \includegraphics[width=.12\linewidth,height=.12\linewidth]{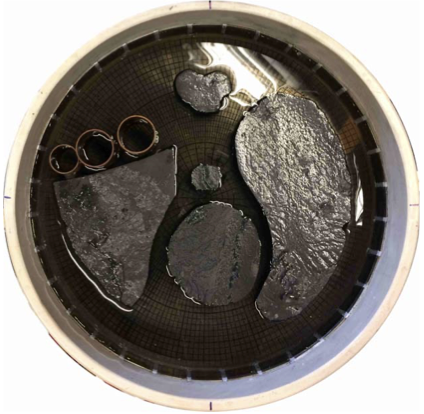}&
    \includegraphics[width=.12\linewidth,height=.12\linewidth]{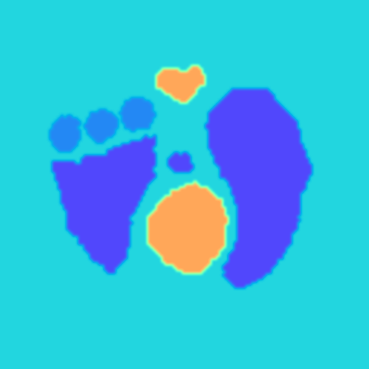}&
    \includegraphics[width=.12\linewidth,height=.12\linewidth]{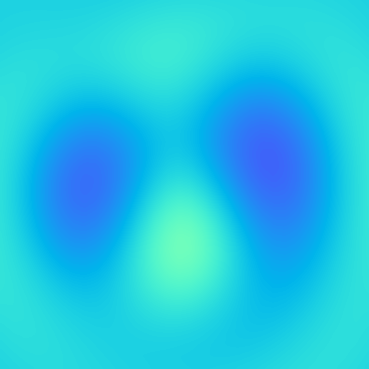}&
    \includegraphics[width=.12\linewidth,height=.12\linewidth]{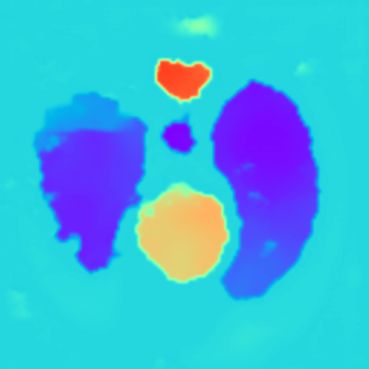}&
    \includegraphics[width=.12\linewidth,height=.12\linewidth]{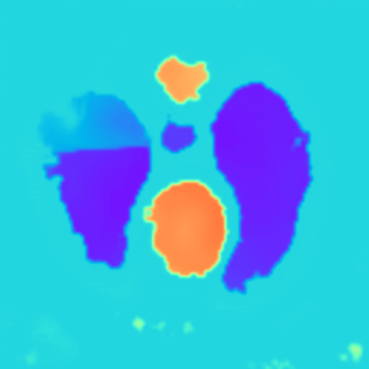}
  \end{tabular}
  }
  \centering
    \begin{minipage}[t]{1\linewidth}
    \hspace{4.4cm} \includegraphics[width=0.5\linewidth]{Figures/colorbar_s.png}
    \caption{In the ACT4 system, the measurements are collected from the complete electrode model. The results are produced from models trained on the simulated data from the continuum model.}
     \label{real_world_ACT4}
    \end{minipage}
\end{figure}

Here, a measurement noise level of $\delta = 0.75\%$ is introduced during the simulation. To ensure the numerical stability and accuracy of the scattering data, we apply a truncation radius $r = 4$ to the scattering data $\mathbf{t}(k)$ for the low-pass regularized D-bar reconstruction. As shown in Figs.~\ref{real_world_KIT4} and \ref{real_world_ACT4}, our cascaded network trained on simulated noisy data demonstrates robust and accurate reconstruction performance for more realistic data.  This contrasts with the single network approach utilized in the Deep D-bar method, which exhibits comparatively less stable and precise results.

\section{\bf{Ablation Study}}\label{6}
\subsection{Multi-Scale Structure}

This study aims to evaluate the numerical performance gain by introducing the multi-scale structure. Here, the multi-scale structure means that the frequency-enhancement U-Net predicts a sequence of high-contrast D-bar images $ \{\sigma_{r_i}(z)\}_{i = 1}^{I}$ for different truncation radii $ \{r_i\}_{i = 1}^{I}$. In contrast, the frequency-enhancement U-Net with the single-scale structure only predicts one high-contrast D-bar image $\sigma_{r_1}(z)$ for a fixed truncation radius $r_1$.

The whole experiment is conducted on the ACT4 dataset with three simulated noise levels $\delta = 0\%, 0.1\%, 0.75\%$. Here, we choose $ r_1, r_2, r_3 = 8, 7, 6$ for the multi-scale structure, while $ r_1 = 7 $ for the single-scale structure. It can be seen from Table~\ref{abla_quan_1} that the multi-scale structure for different measurement noise levels has uniformly lower reconstruction errors compared to the single-scale structure, which demonstrates the effectiveness of the introduced multi-scale structure in the cascaded learning framework.

\begin{table}[htbp]
\centering
\scalebox{0.75}{\begin{tabular}{|cc|cc|cc|cc|}
\hline
\multicolumn{2}{|c|}{\multirow{2}{*}{\textbf{Noise Levels}}}          & \multicolumn{2}{c|}{$\delta = 0.75\%$} & \multicolumn{2}{c|}{$\delta = 0.1\%$} & \multicolumn{2}{c|}{$\delta = 0\%$} \\ \cline{3-8} 
\multicolumn{2}{|c|}{}                                & \multicolumn{1}{c|}{Multi-scale}    & Single-scale   & \multicolumn{1}{c|}{Multi-scale}   & Single-scale   & \multicolumn{1}{c|}{Multi-scale} & Single-scale \\ \hline
\multicolumn{1}{|c|}{\multirow{3}{*}{\textbf{Metrics}}} & PSNR $\uparrow$& \multicolumn{1}{c|}{\textbf{23.63}}     &   23.13   & \multicolumn{1}{c|}{\textbf{25.04} }         &   24.45   & \multicolumn{1}{c|}{\textbf{26.61} }       &  25.89  \\ \cline{2-8} 
\multicolumn{1}{|c|}{}                         & SSIM $\uparrow$& \multicolumn{1}{c|}{\textbf{0.8172} }    &  0.7979    & \multicolumn{1}{c|}{\textbf{0.8521} }    &  0.8345    & \multicolumn{1}{c|}{\textbf{0.8657} }    &  0.8534 \\ \cline{2-8} 
\multicolumn{1}{|c|}{}                         & RMSE $\downarrow$& \multicolumn{1}{c|}{\textbf{0.1445}}   &    0.1542  & \multicolumn{1}{c|}{\textbf{0.1223}}   &    0.1314  & 
\multicolumn{1}{c|}{\textbf{0.1049} }  &    0.1156  \\ \hline
\end{tabular}
}

\caption{Comparison of the reconstruction quality for different U-Net structures }
\label{abla_quan_1}
\end{table}

\subsection{Fixed-Point Iteration Steps}\label{Richardson_exp_part}

\begin{figure}[!htp]
  \centering
  \scalebox{1}{
  \begin{tabular}{ccccccccc}
    &\scriptsize{Ground truth}&\scriptsize{D-bar~\cite{mueller2012linear}}&\scriptsize{1-st step}&\scriptsize{5-th step}&\scriptsize{10-th step}
    \\
    \put(-5,5){\rotatebox{90}{\scriptsize{Sample 1}}}&
    \includegraphics[width=.12\linewidth,height=.12\linewidth]{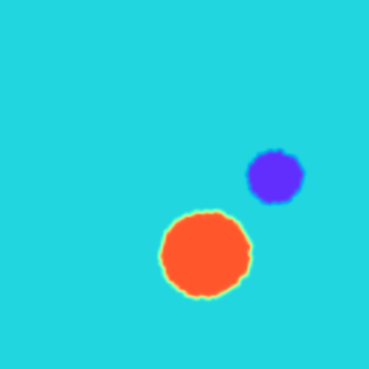}&
    \includegraphics[width=.12\linewidth,height=.12\linewidth]{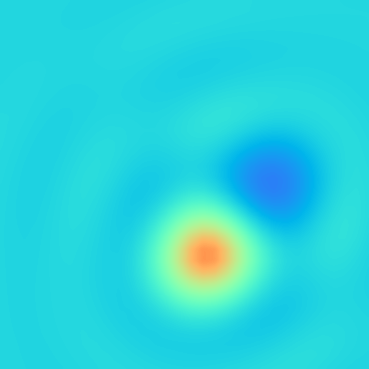}&
    \includegraphics[width=.12\linewidth,height=.12\linewidth]{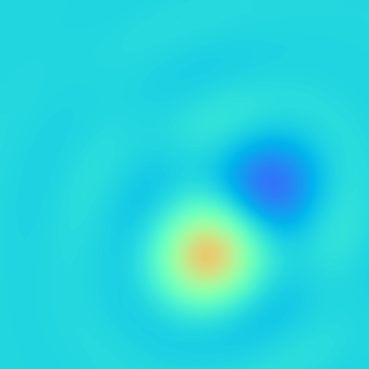}&
    \includegraphics[width=.12\linewidth,height=.12\linewidth]{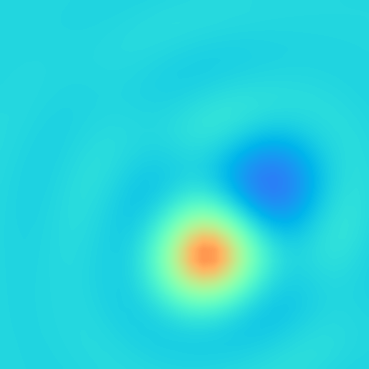}&
    \includegraphics[width=.12\linewidth,height=.12\linewidth]{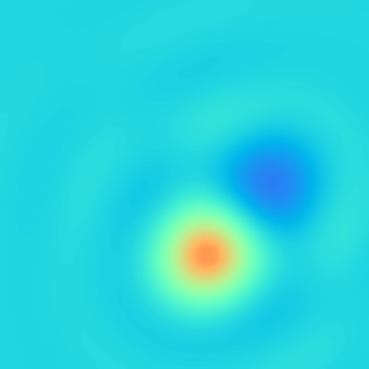}&
    \\ 
    \put(-5,5){\rotatebox{90}{\scriptsize{Sample 2}}}&
    \includegraphics[width=.12\linewidth,height=.12\linewidth]{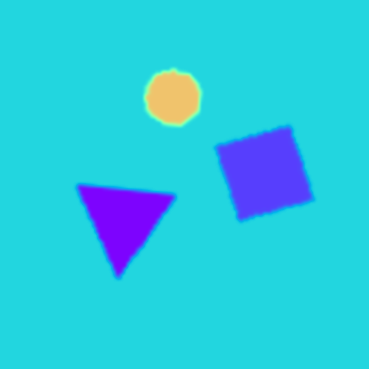}&
    \includegraphics[width=.12\linewidth,height=.12\linewidth]{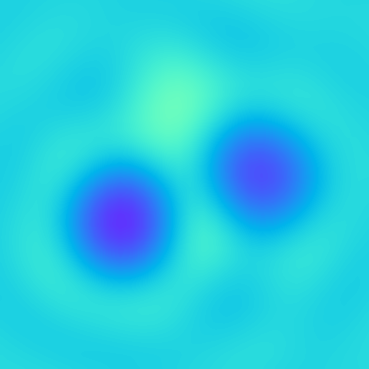}&
    \includegraphics[width=.12\linewidth,height=.12\linewidth]{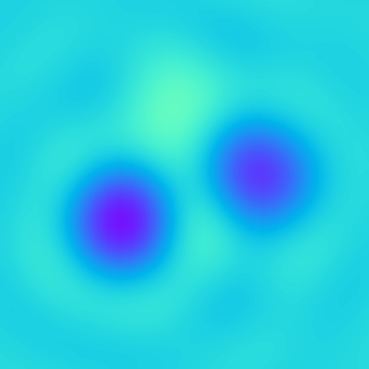}&
    \includegraphics[width=.12\linewidth,height=.12\linewidth]{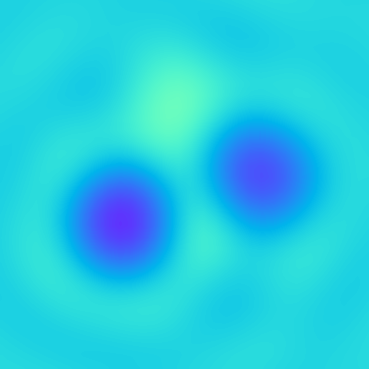}&
    \includegraphics[width=.12\linewidth,height=.12\linewidth]{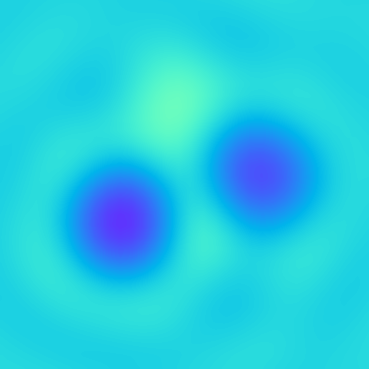}&
    \\ 
    \put(-5,5){\rotatebox{90}{\scriptsize{Sample 3}}}&
    \includegraphics[width=.12\linewidth,height=.12\linewidth]{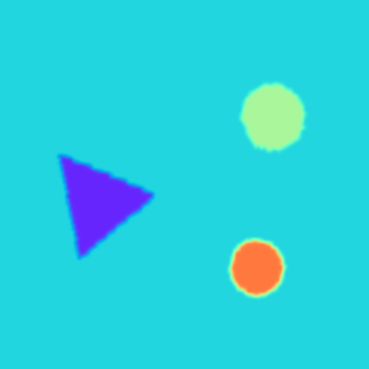}&
    \includegraphics[width=.12\linewidth,height=.12\linewidth]{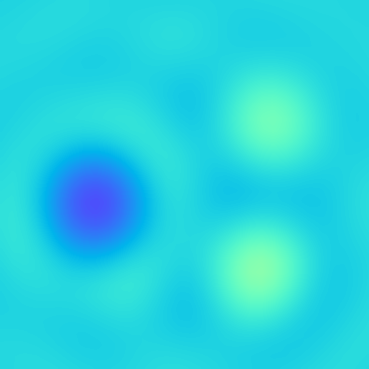}&
    \includegraphics[width=.12\linewidth,height=.12\linewidth]{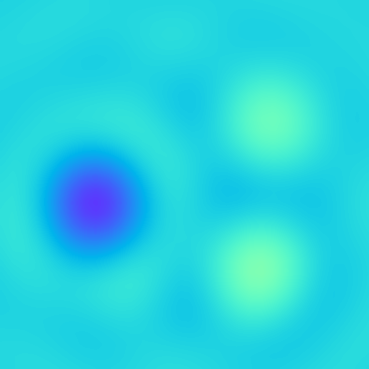}&
    \includegraphics[width=.12\linewidth,height=.12\linewidth]{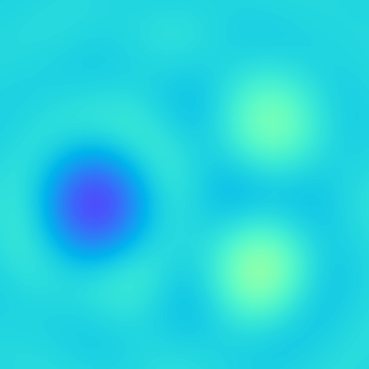}&
    \includegraphics[width=.12\linewidth,height=.12\linewidth]{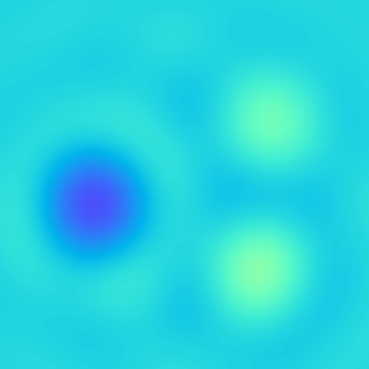}&
    \\ 
    \put(-5,5){\rotatebox{90}{\scriptsize{Sample 4}}}&
    \includegraphics[width=.12\linewidth,height=.12\linewidth]{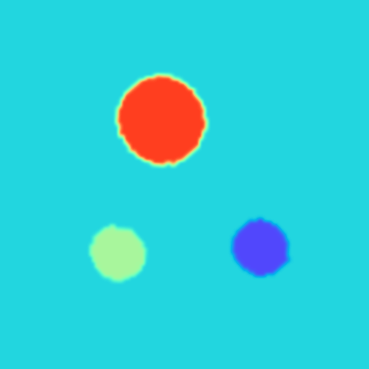}&
    \includegraphics[width=.12\linewidth,height=.12\linewidth]{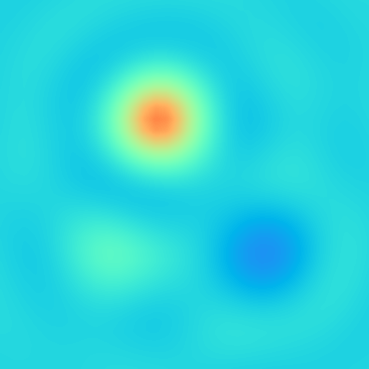}&
    \includegraphics[width=.12\linewidth,height=.12\linewidth]{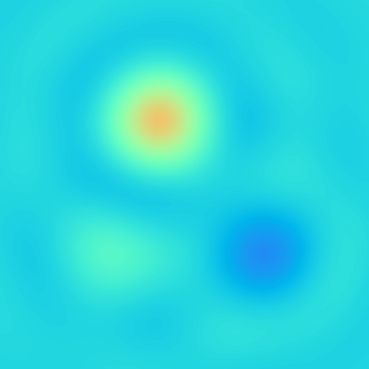}&
    \includegraphics[width=.12\linewidth,height=.12\linewidth]{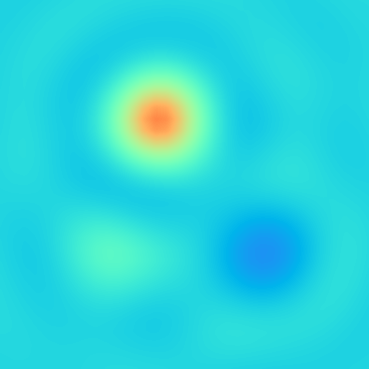}&
    \includegraphics[width=.12\linewidth,height=.12\linewidth]{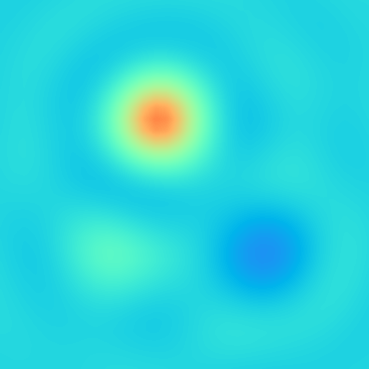}
  \end{tabular}
   }
    \centering
    \begin{minipage}[t]{1\linewidth}
    \hspace{4.4cm} \includegraphics[width=0.5\linewidth]{Figures/colorbar_s.png}
    \caption{The output images are plotted during the fixed-point iteration process for the truncation radius $r = 6$.}
    \label{Richardson_iter}
    \end{minipage}
\end{figure}

We investigate the impact of different iteration steps in the fixed-point method for D-bar reconstruction. Specifically, we examine KIT4 cases with a maximum iteration count of \( S = 1, 5, 10\) in Eq.~(\ref{D-bar-equation-iter}). As illustrated in Fig.~\ref{Richardson_iter}, the results demonstrate that even a single iteration step provides a reasonably accurate approximation of the D-bar reconstruction. A comparison between the reconstructions obtained after the fifth and tenth iterations reveals negligible differences. Based on these observations, we employ a maximum of $5$ iteration steps in our experiments, which proves an effective balance between numerical accuracy and computational efficiency.
    
\section{\bf{Conclusions}}\label{7} 
In this work, we have proposed a multi-scale frequency-enhanced deep D-bar method for real-time Electrical Impedance Tomography (EIT) reconstruction. Building upon the efficiency and simplicity of the regularized D-bar approach, our method addresses its inherent limitations in capturing high-frequency information by integrating a data-driven, deep-learning framework. The proposed methodology combines the GPU-accelerated fixed-point iteration for solving the D-bar integral equations with a cascaded multi-scale U-Net for post-processing, achieving significant improvements in reconstruction quality.

Numerical experiments on simulated KIT4 and ACT4 datasets demonstrate that our method produces high-resolution reconstructions with enhanced contrast and accuracy. 
The GPU-based D-bar acceleration and the direct learning-based post-processing guarantee the computational efficiency for real-time reconstruction. Furthermore, the out-of-distribution testing of the complete electrode model highlights the method's robustness in simulating real-world EIT scenarios.

In summary, our multi-scale frequency-enhanced deep D-bar method offers a novel framework for overcoming the limitations of traditional D-bar reconstructions. By combining physics-informed modeling with data-driven learning, it paves the way for more accurate, efficient, and practical solutions to EIT reconstruction problems.

\clearpage

\bibliographystyle{unsrt}
\bibliography{references}

\newpage
\appendix

\section{Convergence of the Fixed-Point Iteration for D-bar Integral Equations}\label{proof}
\begin{theorem}
Given that  $\sup\limits_{|k| < r} \frac{\left|\mathbf{t}\left(k\right)\right|}{ |k|^2} < \frac{4\pi}{r^2}$, the fixed-point iteration for the D-bar integral equation defined by Eq.~(\ref{D-bar-equation-iter}) converges. Moreover, the error between the $s$-th iterate $m^{(s)}_{r}$ and the exact solution $m_{r}$ satisfies the inequality
$$\sup\limits_{|k|<r} \Big| m^{(s)}_{r}(z,k)  - m_{r}(z,k)\Big| \le \frac{ \lambda^{s+1}_r}{ 1 - \lambda_r},$$
where $\lambda_r = \frac{r^2}{4\pi} \sup\limits_{|k| < r} \frac{\left|\mathbf{t}\left(k\right)\right|}{ |k|^2} $.
\end{theorem}

% For any fixed $z \in \mathbb{R}^2$, the solution of the D-bar integral equation is computed as the solution of
% \begin{equation}\label{A.1}\hspace{0cm}
%     \tilde{m}_R(z, k)= 1 + \frac{1}{(2 \pi)^{2}} \int_{|k^{\prime}|<R} \frac{\tilde{\mathbf{t}}_R\left(k^{\prime}\right)}{\left(k-k^{\prime}\right) \bar{k}^{\prime}} e^{-i(k^{\prime}z + \bar{k}^{\prime}\bar{z})} \overline{\tilde{m}_R\left(z, k^{\prime}\right)} d k^{\prime}.
% \end{equation}

% According to \cite{nachman1996global}, the equation (\ref{A.1}) is a Fredholm equation of the second kind and has a unique solution $\tilde{m}_R$ for any fixed $z \in \mathbb{R}^2$.  

\begin{proof}
To establish the convergence of Eq.~(\ref{D-bar-equation-iter}), we will show that the iterative operator defines a contraction mapping under the specified conditions. Consider the difference between the $s$-th iterate and $(s+1)$-th iterate satisfying
\begin{equation}
     m^{(s+1)}_r(z, k) - m^{(s)}_r(z, k) = \left[\frac{1}{k}\right] \mathbin{\ast} \left( \mathcal{T}_r[m_r^{(s)}] - \mathcal{T}_R[m_r^{(s-1)}] \right).
\end{equation}
Taking the absolute norm $|\cdot |$ on both sides, we expand the explicit formulation in Eq.~(\ref{D-bar-equation}) to obtain 
\begin{equation}\hspace{-0.5cm}
    	 | \Delta m^{(s+1)}_{r}  | = \frac{1}{(2 \pi)^{2}} \left | \int_{|k^{\prime}|<r} \frac{\mathbf{t}\left(k^{\prime}\right)}{\left(k-k^{\prime}\right) \bar{k}^{\prime}} e^{-i(k^{\prime}z + \bar{k}^{\prime}\bar{z})} \overline{\Delta m^{(s)}_r(z,k^\prime)} d k^{\prime}\right |, 
\end{equation}
where we denote $\Delta m^{(s+1)}_r (z, k):= m^{(s+1)}_r(z, k) - m^{(s)}_r(z, k)$, and then
\begin{equation}\hspace{-1.25cm}
    	 | \Delta m^{(s+1)}_{r} | \le \frac{1}{(2 \pi)^{2}}  \int_{|k^{\prime}|<r} \left |\frac{k^{\prime}}{\left(k-k^{\prime}\right)} \right |  \frac{\left |\mathbf{t}\left(k^{\prime}\right)\right |}{|k^{\prime}|^2} \left | e^{-i(k^{\prime}z + \bar{k}^{\prime}\bar{z})}  \right | \left |\overline{\Delta m^{(s)}_r(z,k^\prime)} \right | d k^{\prime}.
\end{equation}
According to Theorem 3.1 in \cite{silt2000dbar}, we know that $|\mathbf{t}(k^{\prime})|\le C |k^{\prime}|^2$ for $k^{\prime}$ near zero. Therefore,  $\frac{\left |\mathbf{t}\left(k^{\prime}\right)\right |}{|k^{\prime}|^2} $ is bounded for $|k^{\prime}| < r$. Continue to take the supremum of $\frac{\left |\mathbf{t}\left(k^{\prime}\right)\right |}{|k^{\prime}|^2} $ and $\left |\overline{\Delta m^{(s)}_r(z,k^\prime)} \right |$ with respect to $k^{\prime}$, and apply $\left |e^{-i(k^{\prime}z + \bar{k}^{\prime}\bar{z})} \right | = 1$ to get
\begin{equation}\label{A.4}\hspace{-1.0cm}
    	 | \Delta m^{(s+1)}_{r}| \le \frac{1}{(2 \pi)^{2}} \sup_{|k^{\prime}| < r} \frac{\left |\mathbf{t}\left(k^{\prime}\right)\right |}{|k^{\prime}|^2} \sup_{|k^{\prime}| < r} \left| \Delta m^{(s)}_r (z, k^\prime)\right | \int_{|k^{\prime}|<r} \frac{|k^{\prime}|}{\left |k-k^{\prime}\right | }  d k^{\prime}.
\end{equation}
Let $C_r:= \sup\limits_{|k^{\prime}| < r} \frac{\left |\mathbf{t}\left(k^{\prime}\right)\right |}{|k^{\prime}|^2}$, and we take the supremum on both sides of Eq.~(\ref{A.4}) for $|k|<r$ to get
\begin{equation}\hspace{0.5cm}
    	 \| \Delta m^{(s+1)}_{r} \|_{\infty} \le \frac{C_r}{(2\pi)^2} \left\| \Delta m^{(s)}_r \right\|_{\infty} \sup\limits_{|k|<r} \int_{|k^{\prime}|<r} \frac{|k^{\prime}|}{\left |k-k^{\prime}\right |}  d k^{\prime}, 
\end{equation}
where $\|\Delta m^{(s+1)}_{r} \|_{\infty}:= \sup\limits_{|k|<r}\Delta m^{(s+1)}_{r}(z,k)$ and $\|\Delta m^{(s)}_{r} \|_{\infty}:= \sup\limits_{|k^{\prime}|<r}\Delta m^{(s)}_{r}(z,k^{\prime})$. Notice that 
\begin{equation}\label{A.6}\hspace{0.75cm}
      \sup\limits_{|k|<r} \int_{|k^{\prime}|<r} \frac{|k^{\prime}|}{\left |k-k^{\prime}\right | }  d k^{\prime} =  \int_{|k^{\prime}|<r} \frac{ \left |k^{\prime}\right |  }{\left | - k^{\prime}\right | }  d k^{\prime} = \pi r^2 ,
\end{equation}
where $\int_{|k^{\prime}|<r} \frac{|k^{\prime}|}{\left |k-k^{\prime}\right | }  d k^{\prime}$ achieves the maximal value at $k = 0$. We can conclude that
\begin{equation}\label{A.7}\hspace{2cm}
  \| \Delta m^{(s+1)}_{r} \|_{\infty} \le \frac{ r^2 C_r }{4 \pi} \left\| \Delta m^{(s)}_r \right\|_{\infty}.
\end{equation}
Therefore, the fixed-point iteration converges provided that:
\begin{equation}\hspace{4cm}
   C_r < \frac{r^2}{4\pi}.
\end{equation}
This condition restricts the choice of the truncation radius $r$ in the regularized D-bar integral equations. By applying the Cauchy Convergence Theorem to Eq.~(\ref{A.7}), we can further conclude that $m^{(s)}_{r}$ will converge to the fixed point $m_{r}$, and the error bound is given by
\begin{equation}\hspace{2.5cm}
    \| m^{(s)}_{r}  - m_{r}  \|_{\infty} \le \frac{ \lambda^{s+1}_r}{ 1 - \lambda_r},
\end{equation}
where $\lambda_r = \frac{r^2}{4\pi} \sup\limits_{|k| < r} \frac{\left|\mathbf{t}\left(k\right)\right|}{ |k|^2}  < 1 $.

\end{proof}

\end{document}